\documentclass[a4paper,11pt,reqno]{amsart}
\usepackage{lmodern}
\usepackage[T1]{fontenc}
\usepackage[utf8]{inputenc}
\usepackage{amsmath}
\usepackage{amsthm}
\usepackage{amssymb}
\usepackage[a4paper]{geometry}
\usepackage{microtype}
\usepackage{xcolor}
\makeatletter

\numberwithin{equation}{section}

\newtheorem{theorem}{Theorem}[section]
\newtheorem{lemma}[theorem]{Lemma}\newtheorem{proposition}[theorem]{Proposition}\newtheorem{corollary}[theorem]{Corollary}

\theoremstyle{definition}
\newtheorem{definition}[theorem]{Definition}

\theoremstyle{remark}
\newtheorem{remark}[theorem]{Remark}

\newcommand{\Hh}{\mathcal H}
\newcommand{\Walg}{\mathcal W}
\newcommand{\CC}{\mathbb C}
\newcommand{\Id}{\operatorname{Id}}
\newcommand{\pr}{\mathrm{pr}}
\newcommand{\Vir}{\mathrm{Vir}}

\title[Whittaker modules  for $\mathcal W^k(\mathfrak{sl}_3)$]{Whittaker modules for the principal $W$-algebra
$\mathcal W^k(\mathfrak{sl}_3)$ via free field realizations}

\author{Dra\v{z}en Adamovi\'{c}}
\email{adamovic@math.hr}
\address{Department of Mathematics, Faculty of Science, University of Zagreb, Bijeni\v{c}ka 30, 10 000 Zagreb, Croatia}

\author{Nina Yu} 
\email{ninayu@xmu.edu.cn}
\address{School of Mathematical Sciences, Xiamen University, Fujian, 361005, China}

\makeatother

\begin{document}
\begin{abstract}
We study Whittaker and generalized Whittaker modules for the
universal principal $W$-algebra $\mathcal W^k(\mathfrak{sl}_3)$,
using two free field realizations: one in the Heisenberg vertex
algebra of rank two and the other in the tensor product of the
universal Virasoro vertex algebra with the Heisenberg vertex
algebra of rank one.
Under explicit conditions on the parameters, we construct universal  Whittaker modules, establish their ordered bases,
and obtain explicit irreducibility criteria.
Under an explicit nondegeneracy condition, we also show that  the submodules of these tensor products correspond to the
submodules of the Virasoro factor. 
We prove irreducibility of the universal Whittaker modules constructed in this paper.

\end{abstract}

\maketitle

\tableofcontents{}

\section{Introduction}\label{introduction}

Whittaker modules have been studied using vertex algebra methods in
several settings. Whittaker modules for the affine Lie algebra
$A_1^{(1)}$ at critical and noncritical levels were constructed
in \cite{ALZ16}. Irreducibility of modules of Whittaker type
for orbifold vertex algebras was studied in \cite{ALPY19,ALPTY25},
including twisted modules and nonabelian orbifolds. Further
constructions concern Whittaker modules for
$\widehat{\mathfrak{gl}}$ and $\mathcal W_{1+\infty}$
\cite{APT23}, and modules of Wakimoto type for
$\widehat{\mathfrak{sl}}_2$ \cite{APT26}.

Whittaker vectors for $W$-algebras have also been studied in connection
with irregular conformal blocks and topological recursion
\cite{KT12,KMST13,BBCC24}. However, we are not aware of a systematic
study of the structure and irreducibility of universal generalized
Whittaker modules for principal $W$-algebras. This is the first in a
series of papers in which we plan to study Whittaker modules for
$W$-algebras.

For $k\in\mathbb C\setminus\{-3\}$, we write
\[
\Walg^k(\mathfrak{sl}_3):=\Walg^k(\mathfrak{sl}_3,f_{\pr})
\]
for the universal principal affine $W$-algebra, where
$f_{\pr}\in\mathfrak{sl}_3$ is a principal nilpotent element.
This is the universal version of Zamolodchikov's algebra \cite{Zam85}.
We write $V^k(\mathfrak{sl}_3)$ for the universal affine vertex
algebra of level $k$. Throughout this paper, we work with
$\Walg^k(\mathfrak{sl}_3)$ rather than its simple quotient,
denoted by $\Walg_k(\mathfrak{sl}_3)$. All modules for vertex
algebras considered below are weak modules.

We use two embeddings of $\Walg^k(\mathfrak{sl}_3)$.
The first takes values in the Heisenberg vertex algebra $\mathcal H_2$
of rank two and is obtained from the type-$A_2$ quantum Miura construction
of Fateev and Lukyanov \cite{FL88}. The construction using two free bosons
also appears in \cite{FZ87}. The second takes values in
$V^{\Vir}(c^{\Vir})\otimes\mathcal H$, where $V^{\Vir}(c^{\Vir})$
is the universal Virasoro vertex algebra and $\mathcal H$ is the
Heisenberg vertex algebra of rank one; it is due to Romans \cite{R91}.
The notation and the homomorphisms $\mu_k$ and $\nu_k$ are defined
in Section~2.

{ One motivation for this work comes from \cite[Remark~7.3]{ALPTY25},  where  it was suggested that certain
 Whittaker modules for the Heisenberg vertex algebras remain irreducible upon restriction
to affine $W$-subalgebras which are subalgebras of  finite orbifolds of Heisenberg vertex algebras. 
At $k=-2$, the Miura embedding identifies
$\Walg^{-2}(\mathfrak{sl}_3)$ with a proper subalgebra of
$\mathcal H_2^{S_3}$, where $S_3$ acts through its standard
two-dimensional representation.
For the modules satisfying the nondegeneracy condition in
Theorem~\ref{thm:rank-two}, irreducibility over
$\mathcal H_2^{S_3}$ follows from
\cite[Theorem~7.1]{ALPTY25}.
We prove that these modules remain irreducible upon restriction
to $\Walg^{-2}(\mathfrak{sl}_3)$ (cf. Remark \ref{orbit-subalgebra}).}

 Although the embedding into $\mathcal H_2$ factors through
$V^{\Vir}(c^{\Vir})\otimes\mathcal H$, we treat the two
constructions separately, {  as each has distinct advantages.}
The direct calculation of the action of $L$ and $W$ on modules for $\mathcal H_2$ provides a starting point for analogous constructions in higher rank.
The realization in $V^{\Vir}(c^{\Vir})\otimes\mathcal H$
allows us to study the submodules in terms of the Virasoro
factor and to obtain precise irreducibility criteria for the
$\Walg^k(\mathfrak{sl}_3)$-modules considered here. 

We put
\[
\alpha_0=\frac{k+2}{\sqrt{k+3}},\qquad
c^{\Vir}=1-6\alpha_0^2.
\]

We say that a pair $(M,w)$ is a Whittaker
$\Walg^k(\mathfrak{sl}_3)$-module of type
$(\chi_1,\chi_2)$ if $w\neq0$ generates $M$ and
\[
L(n)w=\chi_1\delta_{n,1}w,\qquad
W(n)w=\chi_2\delta_{n,1}w\qquad(n\geq1).
\]
For $p\geq1$, we say that $(M,w)$ is a generalized Whittaker
module of order $p$ and type
$(\boldsymbol\ell;\boldsymbol\gamma)
=(\ell_p,\ldots,\ell_{2p};\gamma_{2p},\ldots,\gamma_{3p})$
if $w\neq0$ generates $M$ and
\[
\begin{aligned}
L(n)w&=\ell_nw\quad(p\leq n\leq2p),&
L(n)w&=0\quad(n>2p),\\
W(n)w&=\gamma_nw\quad(2p\leq n\leq3p),&
W(n)w&=0\quad(n>3p).
\end{aligned}
\]
The corresponding universal  Whittaker and generalized Whittaker 
$\Walg^k(\mathfrak{sl}_3)$-modules with these conditions, as defined in
Definition~\ref{def:universal-whittaker}, are denoted by
$\mathbb M^k(\chi_1,\chi_2)$ and
$\mathbb M^k(\boldsymbol\ell;\boldsymbol\gamma)$, respectively.

Under the hypothesis of Theorem~\ref{thm:rank-two}, the homomorphism 
  $\mu_k$ gives irreducible generalized Whittaker
$\Walg^k(\mathfrak{sl}_3)$-modules of arbitrary order.
For the tensor products in Subsection~3.2, we prove that,
under an explicit condition on the parameters, every submodule
after restriction along $\nu_k$ comes from a submodule of the
$\Vir$-module in the first factor. This gives an irreducibility
criterion for those restrictions. For the induced modules used
in these constructions, we also give bases and establish their
universal properties as Whittaker
$\Walg^k(\mathfrak{sl}_3)$-modules.

\begin{theorem}[cf.~Theorem~\ref{thm:whittaker-universal}]\label{intro:whittaker}
Let $\chi_1\neq0$ and $\chi_2\in\mathbb C$.
The  universal Whittaker module $\mathbb M^k(\chi_1,\chi_2)$
exists, has an ordered basis and is irreducible.
 It is realized as a tensor product
of a Virasoro Whittaker module and a highest weight module for the Heisenberg vertex algebra.
\end{theorem}

\begin{theorem}[cf.~Proposition~\ref{prop:generalized-tensor-realization}
and Theorem~\ref{gw:criterion}]\label{intro:generalized}
Let $p\geq1$ and $(\ell_{2p},\gamma_{3p})\neq(0,0)$.
The universal generalized Whittaker module
$\mathbb M^k(\boldsymbol\ell;\boldsymbol\gamma)$ exists, has an
ordered basis and is realized as a tensor product of an induced
generalized Whittaker {module  for the Virasoro
Lie algebra } and an induced module for the Heisenberg
vertex algebra. 
It is irreducible if and only if
\[
\left(\ell_{2p}^{\,3}-\frac{27}{4}\gamma_{3p}^{\,2},\,
2\ell_{2p}\gamma_{3p-1}-3\gamma_{3p}\ell_{2p-1}\right)\neq(0,0).
\]
\end{theorem}

By Corollary~\ref{cor:generalized-tensor-quotients}, every generalized
Whittaker module of this type is obtained by replacing the induced
Virasoro factor by a nonzero quotient. The resulting $\Walg^k(\mathfrak{sl}_3)$-module
is irreducible if and only if that quotient is irreducible.
If $\ell_{2p}^{\,3}-\frac{27}{4}\gamma_{3p}^{\,2}\neq0$,
the universal module is also realized on a module for
$\mathcal H_2$; see Corollary~\ref{cor:rank-two-universal}.

The submodule correspondence also applies when the Virasoro
factor is a highest weight module.

\begin{theorem}[cf.~Theorem~\ref{thm:highest-weight-tensor}]\label{intro:highest-weight}
Let $M$ be a highest weight {module  for the Virasoro Lie algebra} of central charge
$c^{\Vir}$. Let $p\geq1$ and 
$\boldsymbol b=(b_0,\ldots,b_p)$ with $b_p\neq0$, and let $H_p(\boldsymbol b)$ be the induced module for the Heisenberg vertex algebra $\mathcal H$ defined in Section~3.
 Then
$M\otimes H_p(\boldsymbol b)$  is a generalized Whittaker
$\Walg^k(\mathfrak{sl}_3)$-module.
It is irreducible if and only if $M$ is irreducible.
This generalized Whittaker module is  not universal.
\end{theorem}

In particular, for every irreducible highest weight module
$L(c^{\Vir},h)$ for the Virasoro Lie algebra, the tensor product
$L(c^{\Vir},h)\otimes H_p(\boldsymbol b)$ is an irreducible
generalized Whittaker $\Walg^k(\mathfrak{sl}_3)$-module. The first factor
may be a module for a Virasoro minimal model. The tensor product
$M^{\Vir}(c^{\Vir},h)\otimes H_p(\boldsymbol b)$
is not universal, even though its Virasoro factor is a Verma module.

 For convenience, we collect the notation used throughout
the paper and summarize our explicit realizations and their
main properties in Table~\ref{tab:constructions}.
The modules appearing in the table are defined in Section~3.

\subsection*{Notation}

\begin{itemize}

\item
$\Vir$: the Virasoro Lie algebra.

\item $V^{\Vir}(c^{\Vir})$: the universal Virasoro vertex algebra
of central charge $c^{\Vir}$.

\item
$V^k(\mathfrak{sl}_3)$: the universal affine vertex algebra
of level $k$.

\item $\Walg^k(\mathfrak{sl}_3)
=\Walg^k(\mathfrak{sl}_3,f_{\pr})$: the universal principal
affine $W$-algebra. Its simple quotient is denoted by
$\Walg_k(\mathfrak{sl}_3)$.

\item
$\Hh_2$ and $\Hh$: the Heisenberg vertex algebras of ranks
two and one, respectively.

\item
$M^{\Vir}(c^{\Vir},h)$: the Virasoro Verma module of central
charge $c^{\Vir}$ and highest weight $h$.
Its irreducible quotient is denoted by $L(c^{\Vir},h)$.

\item
$M_p^{\Vir}(c^{\Vir},\boldsymbol t)$: the induced generalized
Whittaker $\Vir$-module of order $p$, generated by
$u_{\boldsymbol t}$, where
$\boldsymbol t=(t_p,\ldots,t_{2p})$.

 \item
$L^{\Vir}(c^{\Vir},\chi)
=M_1^{\Vir}(c^{\Vir},(\chi,0))$, $\chi\neq0$:
the irreducible Whittaker module for the Virasoro Lie algebra
with these parameters.

\item
$L^{\Hh_2}(\boldsymbol\lambda)$: the induced $\Hh_2$-module
with parameter
$\boldsymbol\lambda=(a_0,\ldots,a_p;b_0,\ldots,b_p)$.

\item
$H_p(\boldsymbol b)$: the induced $\Hh$-module generated by
$v_{\boldsymbol b}$, where
$\boldsymbol b=(b_0,\ldots,b_p)$.
For $p=0$, we write $F_\beta=H_0(\beta)$.

\end{itemize}

The following table summarizes our explicit realizations and their main properties.

\begin{table}[htbp]
\centering
\caption{Explicit realizations and their properties.}
\label{tab:constructions}
\small
\setlength{\tabcolsep}{4pt}
\renewcommand{\arraystretch}{1.3}

\begin{tabular}{@{}p{0.40\linewidth}p{0.12\linewidth}p{0.12\linewidth}p{0.28\linewidth}@{}}
\hline
Realization  & Condition & Universal & Irreducible \\
\hline

$L^{\mathcal H_2}(\boldsymbol\lambda)$,
restricted along $\mu_k$
&
$\Delta\neq0$
&
Yes
&
Yes
\\[2mm]

$M_p^{\Vir}(c^{\Vir},\boldsymbol t)\otimes H_p(\boldsymbol b)$
&
$D\neq0$
&
Yes
&
If and only if  $ M_p^{\Vir}(c^{\Vir},\boldsymbol t)$  is irreducible
\\[2mm]
$M\otimes H_p(\boldsymbol b)$,
with $M$ a highest weight Vir--module
&
$b_p\neq0$
&
No
&
If and only if $M$ is irreducible
\\[2mm]

$L^{\Vir}(c^{\Vir},\chi_1)\otimes F_\beta$
&
$\chi_1\neq0$
&
Yes
&
Yes
\\
\hline
\end{tabular}

\medskip
\begin{minipage}{\linewidth}
\small
{All properties refer to the resulting
$\Walg^k(\mathfrak{sl}_3)$-modules, obtained by restriction
along $\mu_k$ in the first row and along $\nu_k$ in the
remaining rows. The first three rows concern generalized
Whittaker modules of order $p\geq1$, and the last row concerns
ordinary Whittaker modules. We write
\[
\Delta=\ell_{2p}^{\,3}-\frac{27}{4}\gamma_{3p}^{\,2},
\qquad
D=2t_{2p}-3b_p^2.
\]}
\end{minipage}
\end{table}

The paper is organized as follows. Section~2 recalls $\Walg^k(\mathfrak{sl}_3)$ and the embeddings
$\mu_k:\Walg^k(\mathfrak{sl}_3)\hookrightarrow\mathcal H_2$
and $\nu_k:\Walg^k(\mathfrak{sl}_3)\hookrightarrow
V^{\Vir}(c^{\Vir})\otimes\mathcal H$.
Section~3 introduces the $\mathcal H_2$-modules and
$V^{\Vir}(c^{\Vir})\otimes\mathcal H$-modules used below.
Section~4 defines Whittaker and generalized Whittaker
$\Walg^k(\mathfrak{sl}_3)$-modules.
Section~5 studies restriction along $\mu_k$ and proves irreducibility
under the stated nondegeneracy condition.
Section~6 proves a submodule correspondence for tensor products
of highest weight or generalized Whittaker $\Vir$-modules with
the $\mathcal H$-modules of Section~3, under an explicit condition
on the parameters. It also treats Whittaker
$\Walg^k(\mathfrak{sl}_3)$-modules of type
$(\chi_1,\chi_2)$ for $\chi_1\neq0$.
Section~7 constructs the universal generalized Whittaker modules for
$(\ell_{2p},\gamma_{3p})\neq(0,0)$ through the tensor product
realization, determines their tensor factors, and establishes their
ordered bases.
The proofs of the main results are given in Sections~5--7.

\section{ Free field realizations of  $\Walg^k(\mathfrak{sl}_3)$}
We recall the defining relations of $\Walg^k(\mathfrak{sl}_3)$
and its two free field realizations. We also describe the relation between the corresponding embeddings.
The universal principal affine $W$-algebra is defined by quantum
Drinfeld--Sokolov reduction 
$
\Walg^k(\mathfrak{sl}_3)=H_{\mathrm{DS},f_{\pr}}^{0}\bigl(V^{k}(\mathfrak{sl}_{3})\bigr)$. 
Its central charge is 
$
c_{\Walg}(k)=2-\frac{24(k+2)^{2}}{k+3}$. 
It is strongly generated by fields of weights two and three: 
\[
L(z)=\sum_{n\in\mathbb{Z}}L(n)z^{-n-2},\qquad W(z)=\sum_{n\in\mathbb{Z}}W(n)z^{-n-3}.
\]
Set $c=c_{\Walg}(k)$ and $\Lambda=:L^{2}:-\frac{3}{10}\partial^{2}L$.

We use the $\lambda$-bracket notation as in \cite{Kac98,DSK06}.
With the normalization used below, the defining $\lambda$-brackets are
\begin{align}
[L_\lambda L]&=(\partial+2\lambda)L+\frac{c}{12}\lambda^3,
&[L_\lambda W]&=(\partial+3\lambda)W,\label{eq:LW}\\
[W_\lambda W]&=\frac{c(5c+22)}{17280}\lambda^5
+\frac{5c+22}{144}\lambda^3L
+\frac{5c+22}{96}\lambda^2\partial L\notag\\
&\quad+\lambda\left(\frac23\Lambda+\frac{5c+22}{160}\partial^2L\right)
+\frac13\partial\Lambda+\frac{5c+22}{720}\partial^3L.
\label{eq:WW}
\end{align}
See \cite{A17,Fas25}.
Choose and fix a square root of $k+3$, and set 
$
\alpha_{0}=\frac{k+2}{\sqrt{k+3}}$. 
Thus $c_{\mathcal{W}}\left(k\right)=2-24\alpha_{0}^{2}.$

\subsection{\label{subsec:H2-realization}The embedding $\Walg^k(\mathfrak{sl}_3)\hookrightarrow\mathcal H_2$}

We recall the embedding $\mu_k:\Walg^k(\mathfrak{sl}_3)\hookrightarrow\mathcal H_2$. 
Let $\mathfrak{h}=\mathbb{C}\alpha_{1}\oplus\mathbb{C}\alpha_{2}$
be a two-dimensional vector space with $\langle\alpha_{i},\alpha_{j}\rangle=\delta_{i,j}$
and let

\[
\widehat{\mathfrak{h}}=\mathfrak{h}\otimes\mathbb{C}[t,t^{-1}]\oplus\mathbb{C}K
\]
be the Heisenberg Lie algebra with 
\[
[\alpha_{i}(m),\alpha_{j}(n)]=m\delta_{i,j}\delta_{m+n,0}K,\ [K,\widehat{\mathfrak{h}}]=0,
\]
where $\alpha_{i}(n)=\alpha_{i}\otimes t^{n}.$ Let $\mathcal{H}_{2}$
denote the Heisenberg vertex algebra of rank two associated with
$\widehat{\mathfrak h}$, with $K$ acting as the identity, generated by the fields $\alpha_{i}(z)=\sum_{n\in\mathbb{Z}}\alpha_{i}(n)z^{-n-1},i=1,2.$
The quantum Miura construction of Fateev and Lukyanov gives the
following embedding into $\mathcal H_2$.

\begin{proposition}\label{Fateev=00003D002013Lukyanov free field realization}
There exists an injective vertex algebra homomorphism, called the
Miura map, 
\begin{equation}
\mu_{k}:\Walg^k(\mathfrak{sl}_3)\hookrightarrow\Hh_{2},\qquad k\neq-3,\label{eq:H2-hom}
\end{equation}
determined by 
\begin{equation}
\mu_{k}(L)=\frac{1}{2}:\!\alpha_{1}^{2}\!:+\frac{1}{2}:\!\alpha_{2}^{2}\!:+\frac{\alpha_{0}}{\sqrt{2}}\,\partial\alpha_{1}+\sqrt{\frac{3}{2}}\,\alpha_{0}\,\partial\alpha_{2},\label{eq:H2-L}
\end{equation}

and 
\begin{align}
\mu_{k}(W)= & \frac{1}{\sqrt{6}}:\!\alpha_{1}^{2}\alpha_{2}\!:-\frac{1}{3\sqrt{6}}:\!\alpha_{2}^{3}\!:+\frac{\alpha_{0}}{2}:\!(\partial\alpha_{1})\alpha_{1}\!:+\frac{\alpha_{0}}{\sqrt{3}}:\!(\partial\alpha_{1})\alpha_{2}\!:\nonumber \\
 & -\frac{\alpha_{0}}{2}:\!(\partial\alpha_{2})\alpha_{2}\!:+\frac{\alpha_{0}^{2}}{2\sqrt{2}}\,\partial^{2}\alpha_{1}-\frac{\alpha_{0}^{2}}{2\sqrt{6}}\,\partial^{2}\alpha_{2}.\label{eq:H2-W}
\end{align}

\end{proposition} 
\begin{proof}
The existence of the homomorphism $\mu_k$ follows from the type-$A_{2}$
quantum Miura construction. Its injectivity follows from \cite[Theorem~5.17]{A17}.
Specializing \cite[Corollary~5.20]{A17} to type $A_{2}$ and expressing
the result with respect to the orthonormal basis $\alpha_{1},\alpha_{2}$
gives \eqref{eq:H2-L}--\eqref{eq:H2-W}.
\end{proof}

\begin{remark}  { At \(k=-2\), we have \(\alpha_0=0\) and \(c=2\). The Weyl group \(S_3\) acts on \(\mathcal H_2\) through its standard two-dimensional representation, and the fields \(\mu_{-2}(L)\) and \(\mu_{-2}(W)\) are invariant under this action. Thus
$
\mu_{-2}\bigl(\mathcal W^{-2}(\mathfrak{sl}_3)\bigr)
\subsetneq \mathcal H_2^{S_3}$. 
The inclusion is proper already in conformal weight \(4\), where the two algebras have dimensions \(3\) and \(4\), respectively.

 }
\end{remark}
\subsection{The embedding $\Walg^k(\mathfrak{sl}_3)\hookrightarrow V^{\Vir}(c^{\Vir})\otimes\mathcal H$}

We recall the embedding
$\nu_k:\Walg^k(\mathfrak{sl}_3)\hookrightarrow
V^{\Vir}(c^{\Vir})\otimes\mathcal H$.
Let $V^{\Vir}(c^{\Vir})$ be the universal Virasoro vertex algebra
of central charge $c^{\Vir}=1-6\alpha_{0}^{2},$ with generating field
$T(z)=\sum_{n\in\mathbb{Z}}T(n)z^{-n-2}$. Let $\mathcal{H}$ be
the rank-one Heisenberg vertex algebra generated by $J(z)=\sum_{n\in\mathbb{Z}}J(n)z^{-n-1}$
with $\left[J(m),J(n)\right]=m\delta_{m+n,0}\Id.$

\begin{proposition}\label{prop:vir-heis-hom} There exists a unique
injective vertex algebra homomorphism 
\begin{equation}
\nu_{k}:\Walg^k(\mathfrak{sl}_3)\longrightarrow V^{\Vir}(c^{\Vir})\otimes\Hh\label{eq:VirH-hom}
\end{equation}
determined by 
\begin{equation}
\nu_{k}(L)=T+\frac{1}{2}:\!J^{2}\!:+\sqrt{\frac{3}{2}}\,\alpha_{0}\,\partial J\label{eq:VirH-L}
\end{equation}
and 
\begin{equation}
\nu_{k}(W)=\sqrt{\frac{2}{3}}:\!JT\!:+\frac{\alpha_{0}}{2}\partial T-\frac{1}{3\sqrt{6}}:\!J^{3}\!:-\frac{\alpha_{0}}{2}:\!(\partial J)J\!:-\frac{\alpha_{0}^{2}}{2\sqrt{6}}\partial^{2}J.\label{eq:VirH-W}
\end{equation}
\end{proposition} 
\begin{proof}
Let $\widetilde{L}$ and $\widetilde{W}$ denote the right-hand sides
of \eqref{eq:VirH-L} and \eqref{eq:VirH-W}, respectively.

The Feigin--Fuchs homomorphism \cite{FF90}, \cite[Example~5.22]{A17},
defines the injective vertex algebra homomorphism 
\[
\rho_{\alpha_{0}}:V^{\Vir}(c^{\Vir})\otimes\mathcal{H}\hookrightarrow\mathcal{H}_{2}
\]
determined by 
\[
\rho_{\alpha_{0}}(T)=\frac{1}{2}:\!\alpha_{1}^{2}\!:+\frac{\alpha_{0}}{\sqrt{2}}\partial\alpha_{1},\qquad\rho_{\alpha_{0}}(J)=\alpha_{2}.
\]
A direct substitution gives 
\[
\rho_{\alpha_{0}}(\widetilde{L})=\mu_{k}(L),\qquad\rho_{\alpha_{0}}(\widetilde{W})=\mu_{k}(W).
\]
Since $L$ and $W$ strongly generate $\Walg^k(\mathfrak{sl}_3)$,
we have 
\[
\mu_{k}\bigl(\Walg^k(\mathfrak{sl}_3)\bigr)\subseteq\rho_{\alpha_{0}}\bigl(V^{\Vir}(c^{\Vir})\otimes\mathcal{H}\bigr).
\]
Therefore 
\[
\nu_{k}=\rho_{\alpha_{0}}^{-1}\circ\mu_{k}
\]
is a well-defined vertex algebra homomorphism satisfying \eqref{eq:VirH-L}--\eqref{eq:VirH-W}.
It is unique because $L$ and $W$ strongly generate $\Walg^k(\mathfrak{sl}_3)$,
and it is injective because $\mu_{k}$ is injective. 
\end{proof}
\begin{remark} The homomorphisms satisfy
$\rho_{\alpha_0}\circ\nu_k=\mu_k$. Thus the embedding $\mu_k$ factors
through the vertex algebra $V^{\Vir}(c^{\Vir})\otimes\mathcal H$. 

\end{remark}

\section{\label{free-field-modules}Modules for $\mathcal H_2$ and $V^{\Vir}(c^{\Vir})\otimes\mathcal H$}

 We recall the  modules for the Heisenberg and Virasoro vertex algebras used in the
constructions of Whittaker and generalized Whittaker
$\Walg^k(\mathfrak{sl}_3)$-modules in the following sections.

\subsection{Modules for the Heisenberg vertex algebra $\mathcal H_2$}

\label{subsec:Rank-two-Heisenberg-modules}

\label{rank two heisenberg module}

Fix $p\in\mathbb{Z}_{\ge1}$ and scalars $a_j,b_j\in\CC$ for $0\leq j\leq p$.
Set 
\begin{equation}
\boldsymbol{\lambda}=(a_{0},a_{1},\ldots,a_{p};b_{0},b_{1},\ldots,b_{p}).\label{eq: lambda}
\end{equation}

Let $\widehat{\mathfrak{h}}_{\geq0}=\mathfrak{h}\otimes\mathbb{C}[t]\oplus\mathbb{C}K.$
Define the one-dimensional $\widehat{\mathfrak{h}}_{\geq0}$-module
$\mathbb{C}_{\boldsymbol{\lambda}}=\mathbb{C}e^{\boldsymbol{\lambda}}$
by 
\begin{gather}
\alpha_{1}(j)e^{\boldsymbol{\lambda}}=a_{j}e^{\boldsymbol{\lambda}}\ (0\leq j\leq p),\qquad\alpha_{1}(j)e^{\boldsymbol{\lambda}}=0\text{\ }(j>p);\nonumber \\
\alpha_{2}(j)e^{\boldsymbol{\lambda}}=b_{j}e^{\boldsymbol{\lambda}}\ (0\leq j\leq p),\qquad\alpha_{2}(j)e^{\boldsymbol{\lambda}}=0\text{\ }(j>p);\qquad Ke^{\boldsymbol{\lambda}}=e^{\boldsymbol{\lambda}}.\label{definition relations}
\end{gather}

Define 
\[
L^{\mathcal{H}_{2}}(\boldsymbol{\lambda})=U(\widehat{\mathfrak{h}})\otimes_{U(\widehat{\mathfrak{h}}_{\geq0})}\mathbb{C}_{\boldsymbol{\lambda}}.
\]
Then $L^{\mathcal{H}_{2}}(\boldsymbol{\lambda})$ is a restricted
$\widehat{\mathfrak{h}}$-module of level one, and hence a weak $\mathcal{H}_{2}$-module
generated by $w=w_{\boldsymbol{\lambda}}=1\otimes e^{\boldsymbol{\lambda}};$see
\cite[Section~4.2]{HY19}. The $\mathcal H_2$-module
$L^{\mathcal H_2}(\boldsymbol\lambda)$ is irreducible by the polynomial
action of $\widehat{\mathfrak h}$ described in \cite[Section~5.1]{HY19}.

Fix the lexicographic order on the pairs $(n,i)$, where $n\geq1$
and $i\in\{1,2\}$, by defining $(n,i)>(m,j)$ if either $n>m$, or
$n=m$ and $i>j$. Then $L^{\mathcal{H}_{2}}(\boldsymbol{\lambda})$ has a basis consisting of 
\begin{equation}
\alpha_{i_{1}}(-n_{1})\cdots\alpha_{i_{s}}(-n_{s})w,\qquad(n_{1},i_{1})\geq\cdots\geq(n_{s},i_{s}).\label{eq:H2-PBW}
\end{equation}

Define a linear operator $D_{\Hh_2}$ on $L^{\mathcal H_2}(\boldsymbol\lambda)$ by
\[
D_{\Hh_2}\bigl(\alpha_{i_1}(-n_1)\cdots\alpha_{i_s}(-n_s)w\bigr)
=(n_1+\cdots+n_s)\alpha_{i_1}(-n_1)\cdots\alpha_{i_s}(-n_s)w.
\]
In particular, $D_{\Hh_2}w=0$ and
$[D_{\Hh_2},\alpha_i(-n)]=n\alpha_i(-n)$ for $i=1,2$ and $n\geq1$.
For $N\geq0$, set
\[
E_N^{\Hh_2}=\ker(D_{\Hh_2}-N\Id),\qquad
B_N=\bigoplus_{j=0}^N E_j^{\Hh_2},
\]
and put $B_N=0$ for $N<0$. Thus
$L^{\mathcal H_2}(\boldsymbol\lambda)=\bigoplus_{N\geq0}E_N^{\Hh_2}$.
We say that a nonzero vector in $E_N^{\Hh_2}$ is homogeneous of degree $N$.
For an arbitrary nonzero vector $v$, we define $\deg v$ to be the least
$N$ such that $v\in B_N$.

\subsection{\label{subsec:Virasoro=002013Heisenberg-module}Modules for $V^{\Vir}(c^{\Vir})\otimes\mathcal H$}

We define the highest weight and Whittaker modules for the two vertex
algebras occurring in the embedding $\nu_k$. Throughout,
$c^{\Vir}=1-6\alpha_0^2$.

Let $\Vir$ be the Virasoro Lie algebra with generators $T(n)$,
$n\in\mathbb Z$, and $C_{\Vir}$, subject to
\[
 [T(m),T(n)]=(m-n)T(m+n)
 +\frac{m^3-m}{12}\delta_{m+n,0}C_{\Vir},
 \qquad [C_{\Vir},\Vir]=0.
\]
For $h\in\mathbb C$, the Verma module $M^{\Vir}(c^{\Vir},h)$
is induced from the one-dimensional module for
$\mathbb CC_{\Vir}\oplus\bigoplus_{n\geq0}\mathbb CT(n)$ on which
$C_{\Vir}$ acts as $c^{\Vir}$, $T(0)$ as $h$, and $T(n)$ as zero
for $n\geq1$. Its generating vector is denoted by $u_h$.
It has basis
\[
 T(-n_1)\cdots T(-n_r)u_h,
 \qquad n_1\geq\cdots\geq n_r\geq1,\quad r\geq0.
\]
We call a nonzero quotient of $M^{\Vir}(c^{\Vir},h)$ a highest
weight $\Vir$-module of highest weight $h$.
We denote its irreducible quotient by $L(c^{\Vir},h)$.

For $p\geq1$, put $\Vir_{\geq p}=\bigoplus_{n\geq p}\mathbb CT(n)$.
Given $\boldsymbol t=(t_p,\ldots,t_{2p})$, let
$\mathbb C_{\boldsymbol t,c^{\Vir}}$ be the one-dimensional
$(\Vir_{\geq p}\oplus\mathbb CC_{\Vir})$-module on which
$T(n)$ acts as $t_n$ for $p\leq n\leq2p$, as zero for $n>2p$,
and $C_{\Vir}$ acts as $c^{\Vir}$. 
Define
\[
M_p^{\Vir}(c^{\Vir},\boldsymbol t)
=U(\Vir)\otimes_{U(\Vir_{\geq p}\oplus\mathbb CC_{\Vir})}
\mathbb C_{\boldsymbol t,c^{\Vir}},
\qquad u_{\boldsymbol t}=1\otimes1,
\]
where $U(\Vir)$ denotes the universal enveloping algebra
of the Lie algebra $\Vir$.

We call a nonzero quotient of this induced module a {generalized
Whittaker $\Vir$-module} of type $\boldsymbol t$ and order $p$.
We allow $\boldsymbol t=0$. The induced $\Vir$-module has basis
\begin{equation}\label{eq:Vir-PBW}
 T(p-n_1)\cdots T(p-n_r)u_{\boldsymbol t},
 \qquad n_1\geq\cdots\geq n_r\geq1,\quad r\geq0.
\end{equation}
\begin{theorem}[{\cite[Theorem~7]{LGZ11}}]\label{thm:Vir-induced-simple}
For every $p\geq1$, $c^{\Vir}\in\mathbb C$ and
$\boldsymbol t=(t_p,\ldots,t_{2p})$, the induced module
$M_p^{\Vir}(c^{\Vir},\boldsymbol t)$ is irreducible if and only if
\[
t_{2p}\neq0\quad\text{or}\quad t_{2p-1}\neq0.
\]
\end{theorem}

For $p=1$, the criterion is $(t_1,t_2)\neq(0,0)$.
In particular, $M_1^{\Vir}(c^{\Vir},(0,0))$ is reducible;
its basis includes the vectors $T(0)^j u_{\boldsymbol0}$,
$j\geq0$. The criterion in Theorem~\ref{thm:Vir-induced-simple}
concerns the induced module itself.
The case $p=1$ is the usual induced Whittaker $\Vir$-module. We shall
also use the notation
\[
 L^{\Vir}(c^{\Vir},\chi)
 :=M_1^{\Vir}(c^{\Vir},(\chi,0)),\qquad \chi\neq0,
\]
for this irreducible Whittaker $\Vir$-module; its cyclic vector
is denoted by $u_\chi$.
Every highest weight $\Vir$-module of central charge $c^{\Vir}$
is a quotient of $M_p^{\Vir}(c^{\Vir},0)$ for every $p\geq1$.
All these $\Vir$-modules are restricted and hence are weak modules
for $V^{\Vir}(c^{\Vir})$.

Let $\widehat{\mathfrak a}$ be the Heisenberg Lie algebra generated
by $J(n)$, $n\in\mathbb Z$, and $K$, with
\[
 [J(m),J(n)]=m\delta_{m+n,0}K,\qquad [K,\widehat{\mathfrak a}]=0.
\]
For $p\geq0$ and $\boldsymbol b=(b_0,\ldots,b_p)$, define
$H_p(\boldsymbol b)$ by induction from the one-dimensional module
for $\mathbb CK\oplus\bigoplus_{n\geq0}\mathbb CJ(n)$ generated
by $v_{\boldsymbol b}$, where
\[
 Kv_{\boldsymbol b}=v_{\boldsymbol b},\qquad
 J(n)v_{\boldsymbol b}=b_nv_{\boldsymbol b}\quad(n\geq0),
 \qquad b_n=0\quad(n>p).
\]
This is a weak module for the Heisenberg vertex algebra $\mathcal H$.
It is identified with $\mathbb C[y_1,y_2,\ldots]v_{\boldsymbol b}$ by
\[
 J(-n)=y_n,\qquad J(n)=b_n+n\frac{\partial}{\partial y_n}
 \quad(n\geq1),\qquad J(0)=b_0\Id.
\]
These formulas show that $H_p(\boldsymbol b)$ is irreducible.
For $p=0$ we write $F_\beta=H_0(\beta)$; this is the highest weight
$\mathcal H$-module with $J(0)$-eigenvalue $\beta$.
If some $b_n$ with $n\geq1$ is nonzero, we call $H_p(\boldsymbol b)$
a Whittaker module for the Heisenberg vertex algebra $\mathcal H$. Appending zero entries to
$\boldsymbol b$ does not change this $\mathcal H$-module.

For any highest weight or generalized Whittaker $\Vir$-module $M$
defined above and any $H_p(\boldsymbol b)$, the tensor product
$M\otimes H_p(\boldsymbol b)$ is a weak
$V^{\Vir}(c^{\Vir})\otimes\mathcal H$-module. We regard it as a
$\Walg^k(\mathfrak{sl}_3)$-module by restriction along $\nu_k$.
In particular, for $p\geq1$ put
\[
 V_p(\boldsymbol t,\boldsymbol b)
 =M_p^{\Vir}(c^{\Vir},\boldsymbol t)\otimes H_p(\boldsymbol b),
 \qquad w=u_{\boldsymbol t}\otimes v_{\boldsymbol b}.
\]
The tensor product has basis
\begin{equation}\label{gw:tensor-basis}
 \begin{gathered}
 v_{\boldsymbol n,\boldsymbol m}
 =T(p-n_1)\cdots T(p-n_r)J(-m_1)\cdots J(-m_s)w,\\
 n_1\geq\cdots\geq n_r\geq1,\quad
 m_1\geq\cdots\geq m_s\geq1,\quad r,s\geq0.
 \end{gathered}
\end{equation}
Define the linear operator $D_p$ on this basis by
\[
 D_pv_{\boldsymbol n,\boldsymbol m}
 =\left(\sum_{i=1}^r n_i+\sum_{j=1}^s m_j\right)
 v_{\boldsymbol n,\boldsymbol m}.
\]
Thus $D_pw=0$, $D_pT(p-n)w=nT(p-n)w$ and
$D_pJ(-n)w=nJ(-n)w$ for $n\geq1$. Set
\[
 E_N=\ker(D_p-N\Id),\qquad
 B_N=\bigoplus_{j=0}^N E_j\quad(N\geq0),\qquad B_N=0\quad(N<0).
\]
For $0\neq z\in V_p(\boldsymbol t,\boldsymbol b)$, its degree is
the least $N$ with $z\in B_N$. If $M$ is a quotient of
$M_p^{\Vir}(c^{\Vir},\boldsymbol t)$, we use the image of $B_N$
as the filtration on $M\otimes H_p(\boldsymbol b)$.

The Whittaker modules used below are
\[
 V(\chi,\beta)=L^{\Vir}(c^{\Vir},\chi)\otimes F_\beta
 =V_1((\chi,0),(\beta,0)),\qquad \chi\neq0.
\]
Writing $w=u_\chi\otimes v_\beta$, formulas
\eqref{eq:VirH-L}--\eqref{eq:VirH-W} give
\begin{equation}\label{eq:Whittaker-character}
 L(1)w=\chi w,\qquad
 W(1)w=\chi\left(\sqrt{\frac23}\,\beta-\frac32\alpha_0\right)w,
\end{equation}
and $L(n)w=W(n)w=0$ for $n\geq2$.

\section{Whittaker and generalized Whittaker $\Walg^k(\mathfrak{sl}_3)$-modules}

We now introduce the notions of Whittaker and generalized Whittaker
modules for $\Walg^k(\mathfrak{sl}_3)$ that will be used
in the constructions of the following two sections. 

Following \cite[Definition~3.7]{BBCC24}, we use the following definition.

\begin{definition} \label{def:whittaker} Let $M$ be a weak $\Walg^k(\mathfrak{sl}_3)$-module.
We say that a nonzero vector $v\in M$ is a Whittaker vector of type $(\chi_{1},\chi_{2})\in\CC^{2}$
if 
\[
L(m)v=\chi_{1}\delta_{m,1}v,\qquad W(m)v=\chi_{2}\delta_{m,1}v\qquad(m>0).
\]

We say that a pair $\left(M,v\right)$ is a \emph{Whittaker module of type
$(\chi_{1},\chi_{2})$} if $v$ is a Whittaker vector of this type
and $M$ is generated by $v$ as a $\ensuremath{\Walg^k(\mathfrak{sl}_3)}$-module.
\end{definition}

Motivated by the $\Walg^k(\mathfrak{sl}_3)$ irregular-state conditions in \cite{KMST13}
and by the generalized Whittaker states considered in \cite{KT12},
we introduce the following higher-order version.

\begin{definition}\label{def:generalized-whittaker} Let $p\in\mathbb{Z}_{\ge1}$,
and write 
\[
\boldsymbol{\ell}=(\ell_{p},\ldots,\ell_{2p}),\qquad\boldsymbol{\gamma}=(\gamma_{2p},\ldots,\gamma_{3p}).
\]
Let $M$ be a weak $\Walg^k(\mathfrak{sl}_3)$-module.
We say that a nonzero vector $v\in M$ is a \emph{generalized Whittaker
vector of order $p$ and type $(\boldsymbol{\ell};\boldsymbol{\gamma})$}
if 
\[
L(m)v=\ell_{m}v\quad(p\leq m\leq2p),\qquad L(m)v=0\quad(m>2p),
\]
and 
\[
W(m)v=\gamma_{m}v\quad(2p\leq m\leq3p),\qquad W(m)v=0\quad(m>3p).
\]
No condition is imposed on $L\left(m\right)v$ for $1\le m<p$ or
on $W\left(m\right)v$ for $1\le m<2p.$

We say that a pair $(M,v)$ is a generalized Whittaker module
of order $p$ and type $(\boldsymbol{\ell};\boldsymbol{\gamma})$
if $v$ is such a vector and $M$ is generated by $v$ as a
$\Walg^k(\mathfrak{sl}_3)$-module.
\end{definition}

In particular, a Whittaker vector $v$ of type $(\chi_{1},\chi_{2})$
is a generalized Whittaker vector of order one and type $\boldsymbol{\ell}=(\chi_{1},0),\boldsymbol{\gamma}=(0,0),$
together with the additional condition $W(1)v=\chi_{2}v$.

{ In the next two sections, we construct Whittaker and generalized
Whittaker $\Walg^k(\mathfrak{sl}_3)$-modules using the two
free-field realizations introduced in Section~2. 
In Section~5, we construct generalized Whittaker
$\Walg^k(\mathfrak{sl}_3)$-modules of arbitrary order using
the rank-two Heisenberg realization.
In Section~6, we use the Virasoro--Heisenberg realization
to construct both Whittaker and generalized Whittaker
$\Walg^k(\mathfrak{sl}_3)$-modules.

}

\section{Generalized Whittaker
$\Walg^k(\mathfrak{sl}_3)$-modules from modules for $\mathcal H_2$}

In this section, we construct generalized Whittaker
$\Walg^k(\mathfrak{sl}_3)$-modules by restricting
Whittaker  $\mathcal H_2$-modules along the Miura embedding.
Under an explicit nondegeneracy condition on the parameters,
we prove that the resulting modules are cyclic and irreducible.

Fix $p\in\mathbb{Z}_{\ge1}$ and let $L^{\Hh_{2}}(\boldsymbol{\lambda})$
be the $\Hh_2$-module introduced in Section~\ref{rank two heisenberg module}.
We regard it as a $\Walg^k(\mathfrak{sl}_3)$-module by
restriction along the Miura map $\mu_{k}$ from Proposition~\ref{Fateev=00003D002013Lukyanov free field realization}.
We write $w=w_{\boldsymbol{\lambda}}$ in this section.

By (\ref{eq:H2-L})--(\ref{eq:H2-W}), the corresponding coefficient
operators on $L^{\mathcal{H}_{2}}(\boldsymbol{\lambda})$ are 
\begin{align}
L(m)= & \frac{1}{2}\sum_{r\in\mathbb{Z}}\bigl(:\!\alpha_{1}(r)\alpha_{1}(m-r)\!:+:\!\alpha_{2}(r)\alpha_{2}(m-r)\!:\bigr)\nonumber \\
 & -(m+1)\alpha_{0}\left(\frac{\alpha_{1}(m)}{\sqrt{2}}+\sqrt{\frac{3}{2}}\,\alpha_{2}(m)\right),\label{eq:H2-L-coefficients}
\end{align}
and 
\begin{align}
W(m)= & \frac{1}{\sqrt{6}}\sum_{r+s+t=m}:\!\alpha_{1}(r)\alpha_{1}(s)\alpha_{2}(t)\!:-\frac{1}{3\sqrt{6}}\sum_{r+s+t=m}:\!\alpha_{2}(r)\alpha_{2}(s)\alpha_{2}(t)\!:\nonumber \\
 & -\frac{\alpha_{0}}{2}\sum_{r+s=m}(r+1):\!\alpha_{1}(r)\alpha_{1}(s)\!:-\frac{\alpha_{0}}{\sqrt{3}}\sum_{r+s=m}(r+1):\!\alpha_{1}(r)\alpha_{2}(s)\!:\nonumber \\
 & +\frac{\alpha_{0}}{2}\sum_{r+s=m}(r+1):\!\alpha_{2}(r)\alpha_{2}(s)\!:+\frac{\alpha_{0}^{2}}{2\sqrt{2}}(m+1)(m+2)\alpha_{1}(m)\nonumber \\
 & -\frac{\alpha_{0}^{2}}{2\sqrt{6}}(m+1)(m+2)\alpha_{2}(m).\label{eq:H2-W-coefficients}
\end{align}

A direct calculation using \eqref{eq:H2-L-coefficients}--\eqref{eq:H2-W-coefficients}
and the defining relations of $w$ in (\ref{definition relations})
shows that there exist scalars $\ell_{m}$ and $\gamma_{m}$ such
that 
\[
L(m)w=\ell_{m}w\quad(p\leq m\leq2p),\qquad L(m)w=0\quad(m>2p),
\]
and 
\[
W(m)w=\gamma_{m}w\quad(2p\leq m\leq3p),\qquad W(m)w=0\quad(m>3p).
\]
Thus $w$ is a generalized Whittaker vector of order $p$ and type
$(\boldsymbol{\ell};\boldsymbol{\gamma})$ in the sense of Definition~\ref{def:generalized-whittaker}.

Put $\xi=a_{p}$ and $\eta=b_{p}$ in (\ref{eq: lambda}). Then 
\begin{equation}
\ell_{2p}=\frac{1}{2}(\xi^{2}+\eta^{2}),\qquad\gamma_{3p}=\frac{\eta(3\xi^{2}-\eta^{2})}{3\sqrt{6}}.\label{eq:highest-parameters}
\end{equation}

Put $a_q=b_q=0$ for $q>p$. For $q\ge0$, define 
\begin{equation}
\widehat{\alpha}_{1}(q)=\alpha_{1}(q)-a_{q}\Id,\qquad\widehat{\alpha}_{2}(q)=\alpha_{2}(q)-b_{q}\Id.\label{alphia_i and hat_=00003D00007Balpha_i=00003D00007D}
\end{equation}
Since $\alpha_1(0)$ and $\alpha_2(0)$ act as $a_0\Id$ and $b_0\Id$,
respectively, $\widehat\alpha_i(0)=0$ on
$L^{\mathcal H_2}(\boldsymbol\lambda)$. Set 
\begin{equation}
A\left(\xi,\eta\right)=\begin{pmatrix}\xi & \eta\\[1mm]
\dfrac{\sqrt{6}}{3}\xi\eta & \dfrac{\sqrt{6}}{6}(\xi^{2}-\eta^{2})
\end{pmatrix}.\label{matrix A}
\end{equation}
Then 
\begin{equation}
\det A\left(\xi,\eta\right)=\frac{\sqrt{6}}{6}\xi(\xi^{2}-3\eta^{2}).\label{determinant}
\end{equation}

\begin{lemma}\label{lem:rank-two} Let $N\geq0$ and $n\geq1$.

(i) For every $v\in B_{N}$, we have 
\begin{equation}
L(p-n)v\equiv\bigl(\xi\alpha_{1}(-n)+\eta\alpha_{2}(-n)\bigr)v\pmod{B_{N+n-1}},\label{eq:rank-two-L-positive}
\end{equation}
and 
\begin{equation}
W(2p-n)v\equiv\left(\frac{\sqrt{6}}{3}\xi\eta\,\alpha_{1}(-n)+\frac{\sqrt{6}}{6}(\xi^{2}-\eta^{2})\alpha_{2}(-n)\right)v\pmod{B_{N+n-1}}.\label{eq:rank-two-W-positive}
\end{equation}

(ii) Put 
\[
\mathcal{L}_{n}=L(n+p)-\ell_{n+p}\Id,\qquad\mathcal{D}_{n}=W(n+2p)-\gamma_{n+2p}\Id,
\]
where $\ell_{m}=0$ for $m>2p$ and $\gamma_{m}=0$ for $m>3p$. Then
\begin{equation}
\mathcal{L}_{n}B_{N}\subseteq B_{N-n},\qquad\mathcal{D}_{n}B_{N}\subseteq B_{N-n}.\label{eq:rank-two-lowering}
\end{equation}
Moreover, 
\begin{equation}
\mathcal{L}_{n}v\equiv\bigl(\xi\widehat{\alpha}_{1}(n)+\eta\widehat{\alpha}_{2}(n)\bigr)v\pmod{B_{N-n-1}},\label{eq:rank-two-L-negative}
\end{equation}
and 
\begin{equation}
\mathcal{D}_{n}v\equiv\left(\frac{\sqrt{6}}{3}\xi\eta\,\widehat{\alpha}_{1}(n)+\frac{\sqrt{6}}{6}(\xi^{2}-\eta^{2})\widehat{\alpha}_{2}(n)\right)v\pmod{B_{N-n-1}}.\label{eq:rank-two-W-negative}
\end{equation}
\end{lemma} 
\begin{proof}
Let $w=w_{\boldsymbol{\lambda}}$ as in Subsection \ref{subsec:Rank-two-Heisenberg-modules}.
Set $a_{q}=b_{q}=0$ for $q\in\mathbb Z\setminus\{0,1,\ldots,p\}$. For $q<0$, put $\widehat\alpha_i(q)=\alpha_i(q)$.
Together with (\ref{alphia_i and hat_=00003D00007Balpha_i=00003D00007D}),
we have

\begin{equation}
\alpha_{1}\left(q\right)=a_{q}\Id+\widehat{\alpha}_{1}\left(q\right),\qquad\alpha_{2}\left(q\right)=b_{q}\Id+\widehat{\alpha}_{2}\left(q\right)\ \ \left(q\in\mathbb{Z}\right).\label{alpha_1=00003D00003Dscalar+hatted part}
\end{equation}
For $q\ge0$ and $r\ge1$, we have 
\[
\widehat{\alpha}_{i}(q)w=0,\qquad[\widehat{\alpha}_{i}(q),\alpha_{j}(-r)]=q\delta_{i,j}\delta_{q,r}\Id.
\]
Hence 
\begin{equation}
\widehat{\alpha}_{i}(q)B_{N}\subseteq B_{N-q}\ (q\in\mathbb{Z}).\label{eq:hat-action}
\end{equation}

(i) From \eqref{eq:H2-L}, 
\begin{eqnarray}
L(p-n) & = & \frac{1}{2}\sum_{r+s=p-n}:\!\alpha_{1}(r)\alpha_{1}(s)\!:+\frac{1}{2}\sum_{r+s=p-n}:\!\alpha_{2}(r)\alpha_{2}(s)\!:\nonumber \\
 &  & -(p-n+1)\alpha_{0}\left(\frac{\alpha_{1}(p-n)}{\sqrt{2}}+\sqrt{\frac{3}{2}}\,\alpha_{2}(p-n)\right).\label{L(p-n)}
\end{eqnarray}

Let $v\in B_N$. Substitution of $\alpha_1(q)=a_q\Id+\widehat\alpha_1(q)$
into the first quadratic sum gives 
\begin{alignat*}{1}
\frac{1}{2}\sum_{r+s=p-n}:\!\alpha_{1}(r)\alpha_{1}(s)\!:v & =\frac{1}{2}\sum_{r+s=p-n}a_{r}a_{s}\,v+\sum_{q=0}^{p}a_{q}\widehat{\alpha}_{1}(p-n-q)v\\
 & +\frac{1}{2}\sum_{r+s=p-n}:\!\widehat{\alpha}_{1}(r)\widehat{\alpha}_{1}(s)\!:v.
\end{alignat*}

The first sum belongs to $B_{N}\subseteq B_{N+n-1}$. Every term in
the last sum belongs to $B_{N-r-s}=B_{N+n-p}\subseteq B_{N+n-1}.$

For $0\leq q<p$, the middle sum satisfies $a_{q}\widehat{\alpha}_{1}(p-n-q)v\in B_{N+n+q-p}\subseteq B_{N+n-1}.$
For $q=p$, it is $a_{p}\widehat{\alpha}_{1}(-n)v=\xi\alpha_{1}(-n)v.$
Therefore, 
\[
\frac{1}{2}\sum_{r+s=p-n}:\!\alpha_{1}(r)\alpha_{1}(s)\!:v\equiv\xi\alpha_{1}(-n)v\pmod{B_{N+n-1}}.
\]

Replacing $\alpha_{1},a_{q},\xi$ by $\alpha_{2},b_{q},\eta$ gives
\[
\frac{1}{2}\sum_{r+s=p-n}:\!\alpha_{2}(r)\alpha_{2}(s)\!:v\equiv\eta\alpha_{2}(-n)v\pmod{B_{N+n-1}}.
\]

Finally, for $i=1,2$, $\alpha_{i}(p-n)v\in B_{N}+B_{N+n-p}\subseteq B_{N+n-1}.$
The last line of \eqref{L(p-n)} therefore belongs to $B_{N+n-1}$, and
\[
L(p-n)v\equiv\bigl(\xi\alpha_{1}(-n)+\eta\alpha_{2}(-n)\bigr)v\pmod{B_{N+n-1}}.
\]

For the second congruence in (i), substitution of
(\ref{alpha_1=00003D00003Dscalar+hatted part}) into the first cubic
sum of \eqref{eq:H2-W-coefficients}, applied to $v\in B_N$, gives 
\begin{align}
\sum_{r+s+t=m}:\!\alpha_{1}(r)\alpha_{1}(s)\alpha_{2}(t)\!:v & =\sum_{r+s+t=m}a_{r}a_{s}b_{t}\,v\nonumber \\
 & +\sum_{q,r=0}^{p}\left(2a_{q}b_{r}\widehat{\alpha}_{1}(m-q-r)+a_{q}a_{r}\widehat{\alpha}_{2}(m-q-r)\right)v\nonumber \\
 & +\sum_{q=0}^{p}\sum_{r+s=m-q}\left(2a_{q}:\!\widehat{\alpha}_{1}(r)\widehat{\alpha}_{2}(s)\!:+b_{q}:\!\widehat{\alpha}_{1}(r)\widehat{\alpha}_{1}(s)\!:\right)v\nonumber \\
 & +\sum_{r+s+t=m}:\!\widehat{\alpha}_{1}(r)\widehat{\alpha}_{1}(s)\widehat{\alpha}_{2}(t)\!:v.\label{First cubic sum of W(m)}
\end{align}

Replacing $\alpha_1,a_q$ by $\alpha_2,b_q$ gives the corresponding expansion of
\[
\sum_{r+s+t=m}:\!\alpha_2(r)\alpha_2(s)\alpha_2(t)\!:.
\]

Take $m=2p-n$ in (\ref{First cubic sum of W(m)}). The first line
on the right-hand side belongs to $B_{N}$. By (\ref{eq:hat-action}),
the summands in the third and fourth lines belong to $B_{N+n-p}$
and $B_{N+n-2p}$, respectively. Hence these three lines belong to
$B_{N+n-1}$.

In the second line, $\widehat{\alpha}_{i}(m-q-r)v\in B_{N+n-2p+q+r}.$
Since $0\leq q,r\leq p$, only $q=r=p$ can contribute modulo $B_{N+n-1}$.
Consequently, 
\[
\sum_{r+s+t=2p-n}:\!\alpha_{1}(r)\alpha_{1}(s)\alpha_{2}(t)\!:v\equiv\bigl(2\xi\eta\,\alpha_{1}(-n)+\xi^{2}\alpha_{2}(-n)\bigr)v\pmod{B_{N+n-1}},
\]

and 
\[
\sum_{r+s+t=2p-n}:\!\alpha_{2}(r)\alpha_{2}(s)\alpha_{2}(t)\!:v\equiv3\eta^{2}\alpha_{2}(-n)v\pmod{B_{N+n-1}}.
\]

The quadratic expansion above and \eqref{eq:hat-action} show that
the quadratic and linear expressions in $W(2p-n)$ map $B_N$ into
$B_N+B_{N+n-p}\subseteq B_{N+n-1}$.

Combining these with their coefficients in \eqref{eq:H2-W-coefficients}
gives 
\[
W(2p-n)v\equiv\left(\frac{\sqrt{6}}{3}\xi\eta\,\alpha_{1}(-n)+\frac{\sqrt{6}}{6}(\xi^{2}-\eta^{2})\alpha_{2}(-n)\right)v\pmod{B_{N+n-1}}.
\]

(ii) The same quadratic expansion at $p+n$ gives 
\[
\begin{aligned}\mathcal{L}_{n}v= & \sum_{q=0}^{p}\left(a_{q}\widehat{\alpha}_{1}(p+n-q)+b_{q}\widehat{\alpha}_{2}(p+n-q)\right)v\\
 & +\frac{1}{2}\sum_{r+s=p+n}\left(:\!\widehat{\alpha}_{1}(r)\widehat{\alpha}_{1}(s)\!:+:\!\widehat{\alpha}_{2}(r)\widehat{\alpha}_{2}(s)\!:\right)v\\
 & -(p+n+1)\alpha_{0}\left(\frac{\widehat{\alpha}_{1}(p+n)}{\sqrt{2}}+\sqrt{\frac{3}{2}}\,\widehat{\alpha}_{2}(p+n)\right)v.
\end{aligned}
\]

Here the term multiplying $\Id$ is removed by subtracting $\ell_{p+n}\Id$,
and $a_{p+n}=b_{p+n}=0$. The last two lines belong to $B_{N-p-n}\subseteq B_{N-n-1}$.
In the first line, $p+n-q\geq n$, with equality only for $q=p$. Hence
$\mathcal{L}_{n}B_{N}\subseteq B_{N-n}$ and 
\[
\mathcal{L}_{n}v\equiv\bigl(\xi\widehat{\alpha}_{1}(n)+\eta\widehat{\alpha}_{2}(n)\bigr)v\pmod{B_{N-n-1}}.
\]

Finally, apply the cubic expansions above with $m=2p+n$. The sums
containing only the parameters $a_{q},b_{q}$, with the coefficients
in $W(m)$, give $\gamma_{m}\Id$. Indeed, all other terms annihilate
$w$, since $m-q-r\geq n$ and $m-q\geq p+n$. The quadratic and linear
expressions have no term multiplying $\Id$ because $m>2p$.

After subtracting $\gamma_{2p+n}\Id$, the terms with two or three
operators, as well as the quadratic and linear expressions from $W(2p+n)$,
map $B_{N}$ into $B_{N-p-n}$. Therefore, 
\begin{eqnarray*}
\mathcal{D}_{n}v & \equiv & \frac{1}{\sqrt{6}}\sum_{q,r=0}^{p}\Bigl(2a_{q}b_{r}\widehat{\alpha}_{1}(2p+n-q-r)\\
 &  & +(a_{q}a_{r}-b_{q}b_{r})\widehat{\alpha}_{2}(2p+n-q-r)\Bigr)v\pmod{B_{N-p-n}}.
\end{eqnarray*}

Since $2p+n-q-r\geq n$, every term belongs to $B_{N-n}$. Unless
$q=r=p$, its index is at least $n+1$. Thus $\mathcal{D}_{n}B_{N}\subseteq B_{N-n}$
and 
\[
\mathcal{D}_{n}v\equiv\left(\frac{\sqrt{6}}{3}\xi\eta\,\widehat{\alpha}_{1}(n)+\frac{\sqrt{6}}{6}(\xi^{2}-\eta^{2})\widehat{\alpha}_{2}(n)\right)v\pmod{B_{N-n-1}}.
\]
\end{proof}
The following theorem constructs irreducible generalized Whittaker
$\Walg^k(\mathfrak{sl}_3)$-modules by restricting
$\mathcal H_2$-modules along $\mu_k$, for every type satisfying
the stated nondegeneracy condition.

\begin{theorem}\label{thm:rank-two} Let $p\in\mathbb{Z}_{\ge1}$,
and let 
$
\ell_{p},\ldots,\ell_{2p}, \gamma_{2p},\ldots,\gamma_{3p}
$
be complex numbers such that 
\begin{equation}
\ell_{2p}^{\,3}-\frac{27}{4}\gamma_{3p}^{\,2}\neq0.\label{eq:discriminant}
\end{equation}
Then there exist scalars $a_{j},b_{j}$ $(0\leq j\leq p)$
such that the weak $\Hh_{2}$-module $L^{\Hh_{2}}(\boldsymbol{\lambda})$,
regarded as a $\Walg^k(\mathfrak{sl}_3)$-module through
\eqref{eq:H2-hom}, is generated by $w$ and is irreducible. Moreover,
\[
L(m)w=\ell_{m}w\quad(p\leq m\leq2p),\qquad L(m)w=0\quad(m>2p),
\]
and 
\[
W(m)w=\gamma_{m}w\quad(2p\leq m\leq3p),\qquad W(m)w=0\quad(m>3p).
\]
\end{theorem} 
\begin{proof}
Let $\zeta_{1},\zeta_{2},\zeta_{3}$ be the roots of 
\[
x^{3}-\ell_{2p}x-\gamma_{3p}=0.
\]
Condition \eqref{eq:discriminant} implies that these roots are distinct.
Since $\zeta_{1}+\zeta_{2}+\zeta_{3}=0$, after choosing an ordering
put 
\[
\xi=\frac{\zeta_{1}-\zeta_{2}}{\sqrt{2}},\qquad\eta=\sqrt{\frac{3}{2}}(\zeta_{1}+\zeta_{2})=-\sqrt{\frac{3}{2}}\zeta_{3}.
\]
Then 
\[
\zeta_{1}=\frac{\xi}{\sqrt{2}}+\frac{\eta}{\sqrt{6}},\qquad\zeta_{2}=-\frac{\xi}{\sqrt{2}}+\frac{\eta}{\sqrt{6}},\qquad\zeta_{3}=-\frac{2\eta}{\sqrt{6}}.
\]
Vieta's formulas give 
\[
\ell_{2p}=\frac{1}{2}(\xi^{2}+\eta^{2}),\qquad\gamma_{3p}=\frac{\eta(3\xi^{2}-\eta^{2})}{3\sqrt{6}}.
\]
Moreover, 
\[
\ell_{2p}^{\,3}-\frac{27}{4}\gamma_{3p}^{\,2}=\frac{1}{8}\xi^{2}(\xi^{2}-3\eta^{2})^{2}.
\]
Hence $A(\xi,\eta)$ is invertible.

Set $a_{p}=\xi$ and $b_{p}=\eta$. Let $1\leq r\leq p-1$, and assume
that $a_{p-r+1},\ldots,a_{p}$ and $b_{p-r+1},\ldots,b_{p}$ have
been defined. Collecting in \eqref{eq:H2-L}--\eqref{eq:H2-W} the
summands containing $a_{p-r}$ or $b_{p-r}$ gives 
\[
\begin{pmatrix}\ell_{2p-r}\\[1mm]
\gamma_{3p-r}
\end{pmatrix}=A(\xi,\eta)\begin{pmatrix}a_{p-r}\\[1mm]
b_{p-r}
\end{pmatrix}+\begin{pmatrix}P_{r}\\[1mm]
Q_{r}
\end{pmatrix},
\]
where $P_{r},Q_{r}$ depend only on the scalars already defined. This
equation uniquely determines $a_{p-r},b_{p-r}$. Induction on $r$
determines $a_{p-1},b_{p-1},\ldots,a_{1},b_{1}$. Finally, 
\[
\begin{pmatrix}\ell_{p}\\[1mm]
\gamma_{2p}
\end{pmatrix}=A(\xi,\eta)\begin{pmatrix}a_{0}\\[1mm]
b_{0}
\end{pmatrix}+\begin{pmatrix}P_{p}\\[1mm]
Q_{p}
\end{pmatrix},
\]
where $P_p,Q_p$ depend only on $a_j,b_j$ with
$1\leq j\leq p$, and on $\alpha_0$.
This uniquely determines $a_{0},b_{0}$.

Let $U=\Walg^k(\mathfrak{sl}_3)\cdot w$. We prove by
induction on $N$ that $B_{N}\subset U$. Assume $B_{N-1}\subset U$
and let $v=\alpha_{i}(-n)v_{1}$, where $v_{1}\in B_{N-n}$. By induction,
$v_{1}\in U$. Lemma~\ref{lem:rank-two} gives, modulo $B_{N-1}\subset U$,
\[
L(p-n)v_{1}\equiv\bigl(\xi\alpha_{1}(-n)+\eta\alpha_{2}(-n)\bigr)v_{1},
\]
\[
W(2p-n)v_{1}\equiv\left(\frac{\sqrt{6}}{3}\xi\eta\,\alpha_{1}(-n)+\frac{\sqrt{6}}{6}(\xi^{2}-\eta^{2})\alpha_{2}(-n)\right)v_{1}.
\]
Since $A(\xi,\eta)$ is invertible, both $\alpha_{1}(-n)v_{1}$ and
$\alpha_{2}(-n)v_{1}$ belong to $U$. Thus $U=L^{\Hh_{2}}(\boldsymbol{\lambda})$.

Let $0\neq N\subset L^{\Hh_{2}}(\boldsymbol{\lambda})$ be a $\Walg^k(\mathfrak{sl}_3)$-submodule.
Choose $0\neq v\in N$ of minimal degree $N_{0}$. If $N_{0}=0$,
then $v\in\CC w$, and cyclicity gives $N=L^{\Hh_{2}}(\boldsymbol{\lambda})$.

Assume $N_{0}>0$, and let $v_{N_{0}}$ be the image of $v$ in $B_{N_{0}}/B_{N_{0}-1}$.
Lemma~\ref{lem:rank-two} and minimality give 
\[
\bigl(\xi\widehat{\alpha}_{1}(n)+\eta\widehat{\alpha}_{2}(n)\bigr)v_{N_{0}}=0,
\]
\[
\left(\frac{\sqrt{6}}{3}\xi\eta\,\widehat{\alpha}_{1}(n)+\frac{\sqrt{6}}{6}(\xi^{2}-\eta^{2})\widehat{\alpha}_{2}(n)\right)v_{N_{0}}=0
\]
for every $n\geq1$. Since $A(\xi,\eta)$ is invertible, $\widehat{\alpha}_{1}(n)v_{N_{0}}=\widehat{\alpha}_{2}(n)v_{N_{0}}=0$
for every $n\geq1$.

If a basis vector occurring in $v_{N_{0}}$ contains $\alpha_{i}(-R)$,
choose $R$ maximal and write 
\[
v_{N_{0}}=\sum_{j=0}^{s}\alpha_{i}(-R)^{j}v_{j},\qquad v_{s}\neq0,
\]
where no $v_{j}$ contains $\alpha_{i}(-R)$. Then 
\[
\widehat{\alpha}_{i}(R)v_{N_{0}}=R\sum_{j=1}^{s}j\alpha_{i}(-R)^{j-1}v_{j}\neq0,
\]
contrary to the preceding equality. Hence $N_0=0$, and $L^{\Hh_2}(\boldsymbol\lambda)$ is irreducible as a
$\Walg^k(\mathfrak{sl}_3)$-module. 
\end{proof}
\begin{remark} \label{orbit-subalgebra} {Let $k=-2$.}
For the modules $L^{\mathcal H_2}(\boldsymbol{\lambda})$ considered above, the condition
$
a_p(a_p^2-3b_p^2)\neq0 $
implies that the associated Whittaker function has trivial stabilizer in $S_3$. Hence \cite[Theorem 7.1]{ALPTY25} gives irreducibility over $\mathcal H_2^{S_3}$. Theorem \ref{thm:rank-two}  shows that these modules remain irreducible upon restriction to the proper subalgebra \(\mathcal W^{-2}(\mathfrak{sl}_3)\). This provides an instance of the phenomenon anticipated in  \cite[Remark 7.3]{ALPTY25}.

\end{remark}
\section{Tensor product constructions of $\Walg^k(\mathfrak{sl}_3)$-modules}

In this section, we study $\Walg^k(\mathfrak{sl}_3)$-modules obtained from
tensor products of  modules for the Virasoro and Heisenberg vertex algebras.
Under an explicit condition on the parameters, we describe
their submodules in terms of the Virasoro factor and obtain
an irreducibility criterion.
We also construct irreducible Whittaker modules of type
$(\chi_1,\chi_2)$ with $\chi_1\neq0$.

\subsection{Submodule correspondence under restriction to $\Walg^k(\mathfrak{sl}_3)$}

We prove that, under the condition in the following theorem,
restriction along $\nu_k$ preserves the submodules of the tensor
products defined in Subsection~3.2. In particular, their
irreducibility is determined by the Virasoro factor.
Put $\kappa=\sqrt{2/3}$ and $Q=\sqrt{3/2}\,\alpha_0$.
We extend $t_a$ and $b_a$ by zero outside $p\leq a\leq2p$
and $0\leq a\leq p$, respectively.

\begin{theorem}\label{thm:tensor-submodules}
Let $p\geq1$, let $M$ be a nonzero quotient of
$M_p^{\Vir}(c^{\Vir},\boldsymbol t)$, and suppose that
\[
2t_{2p}-3b_p^2\neq0.
\]
Every $\Walg^k(\mathfrak{sl}_3)$-submodule of
$M\otimes H_p(\boldsymbol b)$, with the action induced by $\nu_k$,
has the form $M'\otimes H_p(\boldsymbol b)$ for a unique
$\Vir$-submodule $M'\subseteq M$.
In particular, $M\otimes H_p(\boldsymbol b)$ is irreducible as a
$\Walg^k(\mathfrak{sl}_3)$-module if and only if $M$ is
irreducible as a $\Vir$-module.
\end{theorem}

\begin{proof}
We first show that a submodule for
$\Walg^k(\mathfrak{sl}_3)$ is invariant under all $T(a)$
and $J(a)$. Use the quotient filtration $B_N$ from
Subsection~3.2. The defining relations give
\begin{equation}\label{gw:bounds}
(T(a)-t_a\Id)B_N\subseteq B_{N+p-a},\qquad
(J(a)-b_a\Id)B_N\subseteq B_{N-a}\qquad(a\in\mathbb Z).
\end{equation}
For the first inclusion, commute $T(a)-t_a\Id$ through the
basis~\eqref{eq:Vir-PBW} and use
\[
 [T(a)-t_a\Id,T(b)]=(a-b)(T(a+b)-t_{a+b}\Id)
 +(a-b)t_{a+b}\Id
 +\frac{a^3-a}{12}c^{\Vir}\delta_{a+b,0}\Id.
\]
The commutator term satisfies the required bound by induction on
the number of factors. A scalar term can be nonzero only when
$a+b\leq2p$, which gives the same bound. Reordering factors with
indices less than $p$ changes their product only in a smaller
filtration space. The second inclusion follows from the polynomial
description of $H_p(\boldsymbol b)$ in Subsection~3.2.
Both inclusions pass to the quotient.
For $r\in\mathbb Z$ and
$S\in\operatorname{End}_{\mathbb C}(M\otimes H_p(\boldsymbol b))$,
write $S\in\mathcal E_r$ if
$SB_N\subseteq B_{N-r}$ for every $N$.  Thus
$\mathcal E_r\mathcal E_s\subseteq\mathcal E_{r+s}$ and
$\mathcal E_{r+1}\subseteq\mathcal E_r$.
Set
\[
A_m=T(p+m)-t_{p+m},\qquad C_m=J(m)-b_m,
\qquad x_m=\binom{A_m}{C_m}\quad(m\in\mathbb Z).
\]
Both entries of $x_m$ belong to $\mathcal E_m$.

Put $t=t_{2p}$ and $d=b_p$.  Expanding
\eqref{eq:VirH-L}--\eqref{eq:VirH-W} around the scalars
$t_a,b_a$ gives
\begin{equation}\label{eq:short-recovery-identity}
y_m=B(t,d)x_m+R_m(x),\qquad
B(t,d)=
\begin{pmatrix}
1&d\\
\kappa d&\kappa t-d^2/\sqrt6
\end{pmatrix},
\end{equation}
where
\[
y_m=\binom{L(p+m)-\ell^0_{p+m}}{W(2p+m)-\gamma^0_{2p+m}}.
\]
Here $\ell^0_q$ and $\gamma^0_q$, for arbitrary integers $q$, denote
the scalar terms obtained by replacing every $T(a)$ by $t_a$
and every $J(a)$ by $b_a$ in those expansions.  The entries of $R_m$
are sums of linear terms with index at least $m+1$, normally
ordered quadratic terms with total index at least $p+m$, and
normally ordered cubic terms with total index $2p+m$.
For example, the first entry is
\[
R_m^L(x)=\sum_{i=0}^{p-1}b_iC_{p+m-i}
-Q(p+m+1)C_{p+m}
+\frac12\sum_{r\in\mathbb Z}:C_rC_{p+m-r}:.
\]
In the expansion of $W(2p+m)$, the terms with index $m$ are precisely
$\kappa dA_m+(\kappa t-d^2/\sqrt6)C_m$; the products
$:JT:$, $:J^3:$, and $:(\partial J)J:$ have the total indices
just described.  Consequently,
\begin{equation}\label{eq:short-recovery-estimate}
R_m(x)\in\mathcal E_{m+1}^{\oplus2}.
\end{equation}
For the substitutions below, all products retain their order from
the expansion: $:U_rV_s:$ means $U_rV_s$ for $r<0$ and
$V_sU_r$ for $r\geq0$, with the cubic product nested in the same
way as $:J(:JJ:):$.
These sums are finite on each vector.  The same is true after
replacing $x_r$ by any pair of operators in $\mathcal E_r$.
On $B_N$, only $h-N\leq r<0$ and $0\leq r\leq N$ can contribute
to $\sum_r:U_rV_{h-r}:$. Applying this observation twice proves
the assertion for the nested cubic sums.

There is one further consequence of the same index count.
If two such families $x$ and $\widetilde x$ agree modulo
$\mathcal E_{r+a}^{\oplus2}$ at every index $r$, then
\begin{equation}\label{eq:short-recovery-improvement}
R_m(x)-R_m(\widetilde x)
\in\mathcal E_{m+a+1}^{\oplus2}\qquad(a\geq0).
\end{equation}
Indeed, in the difference of two products, each summand
contains one factor from the difference of the families;
that factor increases the sum of the indices by $a$.

Let $U$ be a $\Walg^k(\mathfrak{sl}_3)$-submodule. We have
\[
 \det B(t,d)=\frac{2t-3d^2}{\sqrt6}\neq0.
\]
Hence
$B(t,d)^{-1}y_m$ preserves $U$ and agrees with $x_m$ modulo
$\mathcal E_{m+1}^{\oplus2}$.
Equation~\eqref{eq:short-recovery-identity}, together with
\eqref{eq:short-recovery-improvement}, now gives the following
induction: for each $q\geq1$ there is a family $x_m^{[q]}$,
whose entries preserve $U$, such that
\[
x_m^{[q]}\in\mathcal E_m^{\oplus2},\qquad
x_m^{[q]}-x_m\in\mathcal E_{m+q}^{\oplus2}.
\]
The initial family is $B(t,d)^{-1}y_m$; for the inductive step,
replace it by $B(t,d)^{-1}(y_m-R_m(x^{[q]}))$.
For $u\in U\cap B_N$, choosing $m+q>N$ gives
$x_mu=x_m^{[q]}u\in U^{\oplus2}$.
This proves the required invariance under $T(a)$ and $J(a)$.

Write $H_p(\boldsymbol b)=\mathbb C[z_1,z_2,\ldots]v$,
where $J(-i)$ acts by multiplication by $z_i$ and
$(J(i)-b_i)/i$ acts by $\partial/\partial z_i$.
Set
\[
M'=\{a\in M\mid a\otimes v\in U\}.
\]
The invariance just proved implies that $M'$ is a $\Vir$-submodule and $M'\otimes H_p(\boldsymbol b)\subseteq U$.
Conversely, let $\sum_\alpha a_\alpha\otimes z^\alpha v\in U$.
Choose a term of maximal total degree and apply
$\prod_i(\partial/\partial z_i)^{\alpha_i}$.
This gives $\alpha!a_\alpha\otimes v\in U$.
Subtracting $a_\alpha\otimes z^\alpha v$ and repeating shows
that every coefficient belongs to $M'$.
Thus $U=M'\otimes H_p(\boldsymbol b)$, and this expression
also determines $M'$ uniquely.  The converse follows directly
from the formulas for $\nu_k(L)$ and $\nu_k(W)$.
\end{proof}

\subsection{Generalized Whittaker $\Walg^k(\mathfrak{sl}_3)$-modules
from highest weight modules for the Virasoro Lie algebra}\label{subsec:highest-weight-tensor}

\begin{theorem}\label{thm:highest-weight-tensor}

Let $M$ be a highest weight {module  for the Virasoro Lie algebra} of central charge
$c^{\Vir}$, with highest weight vector $u_h$.
Let $p\geq1$ and $b_p\neq0$.
 Then
$
M\otimes H_p(\boldsymbol b)
$
is a generalized Whittaker $\Walg^k(\mathfrak{sl}_3)$-module
of order $p$, generated by $u_h\otimes v_{\boldsymbol b}$.

This module is irreducible if and only if $M$ is irreducible.
This generalized Whittaker module is not universal.
\end{theorem}
\begin{proof}
The Whittaker conditions follow from
\eqref{eq:VirH-L}--\eqref{eq:VirH-W}. The type is given by
\eqref{gw:ell}--\eqref{gw:gamma} with $\boldsymbol t=0$;
in particular,
\[
\ell_{2p}=\frac{b_p^2}{2},\qquad
\gamma_{3p}=-\frac{b_p^3}{3\sqrt6}.
\]
A highest weight $\Vir$-module is a quotient of
$M_p^{\Vir}(c^{\Vir},\boldsymbol 0)$, and the condition in
Theorem~\ref{thm:tensor-submodules} is $-3b_p^2\neq0$.
That theorem gives the irreducibility assertion. It also shows
that the $\Walg^k(\mathfrak{sl}_3)$-submodule generated
by $u_h\otimes v_{\boldsymbol b}$ is $M'\otimes H_p(\boldsymbol b)$
for a $\Vir$-submodule $M'$ containing $u_h$. Thus $M'=M$,
which proves cyclicity.

By Proposition~\ref{prop:generalized-tensor-realization}, the
universal module of this type is
$M_p^{\Vir}(c^{\Vir},\boldsymbol0)\otimes H_p(\boldsymbol b)$.
The canonical surjection onto $M\otimes H_p(\boldsymbol b)$
annihilates
\[
(T(0)-h)u_{\boldsymbol0}\otimes v_{\boldsymbol b}.
\]
This vector is nonzero by the basis \eqref{eq:Vir-PBW}.
Consequently, the quotient is not universal.
\end{proof}

In particular, $L(c^{\Vir},h)\otimes H_p(\boldsymbol b)$ is an
irreducible generalized Whittaker module whenever $b_p\neq0$.
If $c^{\Vir}$ is the central charge of a Virasoro minimal model,
the first factor may be any irreducible module for that minimal
model. These are modules for the universal algebra
$\Walg^k(\mathfrak{sl}_3)$.
The same construction applies to the Verma module
$M^{\Vir}(c^{\Vir},h)$; the resulting generalized Whittaker
module is not universal.

\begin{remark}\label{rem:reducible-tensors}
It is easy to see that when both factors are irreducible highest weight modules, their
tensor product need not remain irreducible after restriction
to $\Walg^k(\mathfrak{sl}_3)$. But the same  can occur when both factors are irreducible
Whittaker modules and the condition in
Theorem~\ref{thm:tensor-submodules} is not satisfied.
Take $k=-2$, so that $\alpha_0=0$ and $c^{\Vir}=1$, and consider
\[
M_1^{\Vir}(1,(0,3/2))\otimes H_1(0,1).
\]
Both factors are irreducible, but $2t_2-3b_1^2=0$.
To see that the restriction is reducible, realize the first factor
on a module for the Heisenberg vertex algebra of rank one,
generated by $u$, with
\[
A(0)u=0,\qquad A(1)u=-\sqrt3\,u,\qquad A(n)u=0\quad(n\geq2),
\qquad T=\tfrac12:A^2:.
\]
Then $T(1)u=0$ and $T(2)u=\tfrac32u$.
For the filtration defined by $\deg A(-n)=n$, the leading
action of $T(1-n)$ is multiplication by $-\sqrt3 A(-n)$.
Thus the basis \eqref{eq:Vir-PBW} gives exactly this module. Write $B=J$ for the second Heisenberg field
and put
\[
X=\frac{A+\sqrt3 B}{2},\qquad
Y=\frac{-\sqrt3 A+B}{2}.
\]
The tensor product becomes $F_0^X\otimes H_1^Y(0,2)$, and
\[
L=\frac12:X^2:+\frac12:Y^2:,
\qquad
W=-\frac1{\sqrt6}:X^2Y:+\frac1{3\sqrt6}:Y^3:.
\]
Both fields are invariant under $X\mapsto-X$, $Y\mapsto Y$.
The even and odd subspaces for this involution are nonzero
$\Walg^k(\mathfrak{sl}_3)$-submodules, proving reducibility.
\end{remark}

\subsection[Generalized Whittaker $\Walg^k(\mathfrak{sl}_3)$-modules from induced Virasoro Whittaker modules]{Generalized Whittaker $\Walg^k(\mathfrak{sl}_3)$-modules\ from induced Virasoro Whittaker modules}

\begin{proposition}\label{prop:general-tensor-basis}
Let $p\geq1$ and $2t_{2p}-3b_p^2\neq0$.
The $\Walg^k(\mathfrak{sl}_3)$-module
$V_p(\boldsymbol t,\boldsymbol b)$ is generated by $w$ and has basis
\[
 \begin{gathered}
 L(p-n_1)\cdots L(p-n_r)W(2p-m_1)\cdots W(2p-m_s)w,\\
 n_1\geq\cdots\geq n_r\geq1,\quad
 m_1\geq\cdots\geq m_s\geq1,\quad r,s\geq0.
 \end{gathered}
\]
\end{proposition}
\begin{proof}
Apply \eqref{eq:short-recovery-identity} with $m=-n$, $n\geq1$.
For $z\in B_N$, the estimate \eqref{eq:short-recovery-estimate} gives
\begin{equation}\label{gw:leading}
 \binom{L(p-n)z}{W(2p-n)z}
 \equiv B(t,d)\binom{T(p-n)z}{J(-n)z}
 \pmod{B_{N+n-1}^{\oplus2}}.
\end{equation}
The scalar terms also belong to $B_N\subseteq B_{N+n-1}$.
By \eqref{gw:tensor-basis}, the associated graded space is
$\mathbb C[x_n,y_n\mid n\geq1]$, where $x_n,y_n$ are the
classes of $T(p-n)$ and $J(-n)$ and have degree $n$.
The Virasoro commutation relations give
\[
 [T(p-n),T(p-r)]B_N\subseteq B_{N+n+r-1}
 \qquad(n,r\geq1).
\]
Thus $T(p-n)$ and $J(-n)$ induce multiplication by $x_n$ and
$y_n$, respectively. Formula~\eqref{gw:leading} identifies the
classes of the displayed $L,W$-vectors with the polynomial basis
after the changes of variables
\[
 \binom{x_n}{y_n}\longmapsto B(t,d)\binom{x_n}{y_n}
 \qquad(n\geq1).
\]
Since $B(t,d)$ is invertible, those classes form a basis of
$B_N/B_{N-1}$ for every $N$. Induction on $N$ proves that the
displayed vectors form a basis of the tensor product and that $w$
generates it as a $\Walg^k(\mathfrak{sl}_3)$-module.
\end{proof}

\begin{corollary}\label{cor:general-tensor-simple}
Suppose that $p\geq1$ and $2t_{2p}-3b_p^2\neq0$.
Then $V_p(\boldsymbol t,\boldsymbol b)$ is an irreducible
$\Walg^k(\mathfrak{sl}_3)$-module if and only if
\[
 (t_{2p},t_{2p-1})\neq(0,0).
\]
\end{corollary}
\begin{proof}
Theorem~\ref{thm:tensor-submodules} reduces irreducibility
to the induced $\Vir$-module. The assertion follows from
Theorem~\ref{thm:Vir-induced-simple}.
\end{proof}

\subsection{Irreducible Whittaker $\Walg^k(\mathfrak{sl}_3)$-modules}

We consider the Whittaker modules
$V(\chi,\beta)=L^{\Vir}(c^{\Vir},\chi)\otimes F_\beta$,
where $\chi\neq0$. Here $t_2=b_1=0$, so the preceding theorem
does not apply. We use the filtration $B_N$ of Subsection~3.2
with $p=1$.

\begin{lemma}\label{lem:vir-heis} Let $N\geq0$ and $n\geq1$, and
set $\kappa=\sqrt{\frac{2}{3}}.$

(i) For every $v\in B_{N}$, we have 
\begin{equation}
L(1-n)v\equiv T(1-n)v\pmod{B_{N+n-1}},\label{eq:vir-heis-L-positive}
\end{equation}
and 
\begin{align}
W(1-n)v\equiv{} & \kappa\chi J(-n)v+\left(\kappa\beta+\frac{n-3}{2}\alpha_{0}\right)T(1-n)v\nonumber \\
 & +\kappa\sum_{j=1}^{n-1}J(-j)T(1-n+j)v\nonumber \\
 & +\kappa\sum_{r\geq1}T(1-n-r)J(r)v\pmod{B_{N+n-1}}.\label{eq:vir-heis-W-positive}
\end{align}
The last sum is finite on every vector.

(ii) Define 
\[
\mathcal{L}_{n}=L(n)-\delta_{n,1}\chi\Id,\qquad\mathcal{D}_{n}=W(n+1).
\]
Then 
\begin{equation}
\mathcal{L}_{n}B_{N}\subseteq B_{N-n},\qquad\mathcal{D}_{n}B_{N}\subseteq B_{N-n}.\label{eq:vir-heis-lowering}
\end{equation}
Moreover, for every $v\in B_{N}$, 
\begin{equation}
\mathcal{D}_{n}v\equiv\kappa\left(\chi J(n)+\sum_{r>n}T(1-r+n)J(r)\right)v\pmod{B_{N-n-1}}.\label{eq:vir-heis-D-negative}
\end{equation}
\end{lemma}  
\begin{proof}
Equation \eqref{eq:VirH-L} gives the first congruence. In \eqref{eq:VirH-W},
the relevant summands of $:\!JT\!:$ are $J(0)T(1-n)$, $J(-n)T(1)$,
$J(-j)T(1-n+j)$, $1\leq j\leq n-1$, and $J(r)T(1-n-r)$, $r\geq1$.
Since $T(1)v\equiv\chi v\pmod{B_{N-1}}$, these summands give the
required congruence. The summand $\frac{\alpha_{0}}{2}\partial T$
contributes $\frac{n-3}{2}\alpha_{0}T(1-n)v$. The remaining summands
in \eqref{eq:VirH-W} belong to $B_{N+n-1}$.

For $\mathcal{D}_{n}$, the relevant summands of $:\!JT\!:$ are $J(n)T(1)$
and $J(r)T(1-r+n)$, $r>n$. They give the last congruence. The remaining
summands belong to $B_{N-n-1}$. The two inclusions follow from \eqref{eq:VirH-L}--\eqref{eq:VirH-W}. 
\end{proof}
\begin{theorem}\label{thm:vir-heis} Let $(\chi_{1},\chi_{2})\in\CC^{2}$
with $\chi_{1}\neq0$, and put $\beta=\sqrt{\frac{3}{2}}\left(\frac{\chi_{2}}{\chi_{1}}+\frac{3}{2}\alpha_{0}\right).$
Then the tensor product 
\[
L^{\Vir}(c^{\Vir},\chi_{1})\otimes F_{\beta}
\]
is an irreducible $\Walg^k(\mathfrak{sl}_3)$-module under \eqref{eq:VirH-hom}.
It is generated by the Whittaker vector $w=u_{\chi_{1}}\otimes v_{\beta}$
satisfying 
\[
L(1)w=\chi_{1}w,\qquad W(1)w=\chi_{2}w,\qquad L(n)w=W(n)w=0,\qquad n\geq2.
\]
\end{theorem} 
\begin{proof}
Equation \eqref{eq:Whittaker-character} gives the stated Whittaker
conditions.

Let $U=\Walg^k(\mathfrak{sl}_3)\cdot w$.
For a tensor basis vector in \eqref{gw:tensor-basis}, let
$(N,s,S)$ denote its degree, the number of $J$-factors,
and the sum of their indices, respectively.
We prove that every such vector belongs to $U$ by induction on
$(N,s,S)$ with respect to the lexicographic order.

If no factor $J(-m)$ occurs, the assertion follows from the first
congruence in Lemma~\ref{lem:vir-heis}. Suppose that $J(-n)v$ is
a vector in \eqref{gw:tensor-basis} and that the assertion is known
for vectors of smaller degree, for vectors of the same degree with
fewer factors $J(-m)$, and, when the number of these factors is
fixed, for a smaller sum of their indices. Put $N=\deg(J(-n)v)$.
By induction hypothesis,
$v\in U$. The second congruence in Lemma~\ref{lem:vir-heis}, modulo
$B_{N-1}\subset U$, expresses $\kappa\chi_{1}J(-n)v$ as a linear
combination of vectors belonging to $U$. Since $\chi_{1}\neq0$,
$J(-n)v\in U$. Hence $U=V(\chi_{1},\beta)$.

Let $0\neq M\subset V(\chi_{1},\beta)$ be a $\Walg^k(\mathfrak{sl}_3)$-submodule.
Choose $0\neq v\in M$ of minimal degree $N_{0}$, and let $\overline{v}$
be its image in $B_{N_{0}}/B_{N_{0}-1}$. Minimality gives $\mathcal{L}_{n}v=\mathcal{D}_{n}v=0$
for every $n\geq1$.

Assume that a factor $J(-R)$ occurs in a basis vector appearing in
$\overline{v}$, with $R$ maximal. Write 
\[
\overline{v}=\sum_{j=0}^{s}J(-R)^{j}v_{j},\qquad v_{s}\neq0,
\]
where no $v_{j}$ contains $J(-R)$. Since $J(r)\overline{v}=0$ for
$r>R$, the last congruence in Lemma~\ref{lem:vir-heis}, with $n=R$,
gives $\chi_{1}J(R)\overline{v}=0$. However, 
\[
J(R)\overline{v}=R\sum_{j=1}^{s}jJ(-R)^{j-1}v_{j}\neq0.
\]
This is a contradiction. Hence $\overline{v}=\widetilde{v}\otimes v_{\beta}$
for some $\widetilde{v}\in L^{\Vir}(c^{\Vir},\chi_{1})$.

Write the leading class as a polynomial in the variables
$x_r=T(1-r)$, $r\geq1$. If $N_0>0$, choose the largest $R$
for which $x_R$ occurs. On this polynomial, the component of
$\mathcal L_R$ in $B_{N_0-R}/B_{N_0-R-1}$ is
\[
 (2R-1)\chi_1\frac{\partial}{\partial x_R}.
\]
Indeed, $[T(R),T(1-R)]=(2R-1)T(1)$, and $T(1)$ contributes
$\chi_1$ in that quotient. A commutator with $T(1-r)$ for $r<R$
produces $T(a)$ with $a\geq2$; commuting it through another factor
gives a vector in $B_{N_0-R-1}$. The terms involving $J(n)$
contribute zero to the leading class, which has no $J(-n)$ factors.
The displayed derivative is nonzero, contradicting
$\mathcal L_Rv=0$. Hence $N_0=0$, and $v$ is a nonzero scalar
multiple of $w$. Cyclicity now gives $M=V(\chi_1,\beta)$.
\end{proof}

\section{Universal $\Walg^k(\mathfrak{sl}_3)$-modules of Whittaker type}\label{sec:universal}

For $(\ell_{2p},\gamma_{3p})\neq(0,0)$, we construct the universal
generalized Whittaker modules, determine their tensor factors,
and express irreducibility in terms of the  Whittaker conditions.

\subsection{Universal Whittaker $\Walg^k(\mathfrak{sl}_3)$-modules and their ordered bases}

\begin{definition}\label{def:universal-whittaker}
Let $(M,w)$ be a Whittaker or generalized Whittaker
$\Walg^k(\mathfrak{sl}_3)$-module.
We say that $(M,w)$ is \emph{universal} if, for every weak module
$X$ containing a vector $v$ with the same Whittaker conditions,
there is a unique module homomorphism $M\to X$ taking $w$ to $v$.
In the generalized case, the order $p$ is fixed.
We denote the universal modules, when they exist, by
$\mathbb M^k(\chi_1,\chi_2)$ and
$\mathbb M^k(\boldsymbol\ell;\boldsymbol\gamma)$, respectively.
\end{definition}

Every module generated by a vector with the prescribed conditions
is a quotient of the corresponding universal module.

\begin{definition}\label{def:PBW-whittaker}
Put
\begin{equation}\label{eq:r_L r_W}
(r_L,r_W)=
\begin{cases}
(1,1),&\text{for Whittaker modules},\\
(p,2p),&\text{for generalized Whittaker modules of order }p.
\end{cases}
\end{equation}
Let $\mathcal S$ consist of the ordered products
\[
L(r_L-n_1)\cdots L(r_L-n_r)
W(r_W-m_1)\cdots W(r_W-m_s),
\]
where $n_1\geq\cdots\geq n_r\geq1$,
$m_1\geq\cdots\geq m_s\geq1$, and $r,s\geq0$.
The empty product is the identity.
We say that $(M,w)$ \emph{admits an ordered basis} if
$(Pw)_{P\in\mathcal S}$ is a basis of $M$.
\end{definition}

\begin{lemma}\label{lem:Whittaker-spanning}
If $x$ is a Whittaker or generalized Whittaker vector, then
the module generated by $x$ is spanned by $Px$, $P\in\mathcal S$,
with the corresponding values of $r_L,r_W$.
\end{lemma}
\begin{proof}
Put $h(L)=2$, $h(W)=3$, and let $F_d$ be the span of the vectors
\[
A_1(n_1)\cdots A_r(n_r)x,\qquad
A_i\in\{L,W\},\quad n_i\in\mathbb Z,\quad
\sum_i h(A_i)\leq d.
\]
Set $F_d=0$ for $d<0$. Strong generation by $L,W$ gives
$M=\bigcup_{d\geq0}F_d$ for the module $M$ generated by $x$.
The modes of $\Lambda=:L^2:-\frac3{10}\partial^2L$ are
locally finite sums of products of at most two $L$-modes.
The defining $\lambda$-brackets therefore imply
\[
[A(m),B(n)]F_d
\subseteq F_{d+h(A)+h(B)-2}
\qquad(A,B\in\{L,W\}).
\]
Consequently, any permutation of the factors of a word of total
weight at most $d$ changes its value only by an element of $F_{d-2}$.
For $n\geq r_A$, the vector $A(n)x$ is a prescribed scalar
multiple of $x$. Ordering the remaining factors gives
\[
F_d=\operatorname{Span}\{Px:P\in\mathcal S,\ h(P)\leq d\}
     +F_{d-2},
\]
where $h(P)$ is the sum of the weights of the factors of $P$.
Induction on $d$ proves the assertion.
\end{proof}

\begin{proposition}\label{prop:PBW-implies-universal}
A Whittaker or generalized Whittaker module admitting the basis
in Definition~\ref{def:PBW-whittaker} is universal.
\end{proposition}
\begin{proof}
Let $(M,w)$ admit this basis, and let $v\in X$ satisfy the same
Whittaker conditions. The submodule $Y\subseteq M\oplus X$
generated by $(w,v)$ satisfies
\[
Y=\operatorname{Span}\{(Pw,Pv):P\in\mathcal S\}
\]
by Lemma~\ref{lem:Whittaker-spanning}.
Since $(Pw)_{P\in\mathcal S}$ is a basis of $M$, the first
projection $Y\to M$ is an isomorphism. Its inverse followed by
the second projection gives a homomorphism $M\to X$ taking
$w$ to $v$. It is unique because $w$ generates $M$.
\end{proof}

\subsection{Factorization of generalized Whittaker $\Walg^k(\mathfrak{sl}_3)$-modules}\label{gw:subsection}

For
\[
V_p(\boldsymbol t,\boldsymbol b)
=M_p^{\Vir}(c^{\Vir},\boldsymbol t)\otimes H_p(\boldsymbol b),
\]
write $w=u_{\boldsymbol t}\otimes v_{\boldsymbol b}$,
$\kappa=\sqrt{2/3}$ and $Q=\sqrt{3/2}\alpha_0$.
Extend $t_a,b_a$ by zero outside $p\leq a\leq2p$
and $0\leq a\leq p$, respectively.
The realization \eqref{eq:VirH-L}--\eqref{eq:VirH-W} gives
\begin{align}
\ell_m&=t_m+\frac12\sum_{i=0}^{p}b_ib_{m-i}
             -\delta_{m,p}Q(p+1)b_p,
             &&p\leq m\leq2p,\label{gw:ell}\\
\gamma_m&=\kappa\sum_{i=0}^{p}b_it_{m-i}
-\frac1{3\sqrt6}\sum_{i,j=0}^{p}b_ib_jb_{m-i-j}
+\delta_{m,2p}\alpha_0(p+1)
       \left(\frac{b_p^2}{2}-t_{2p}\right),
             &&2p\leq m\leq3p.\label{gw:gamma}
\end{align}
Thus $w$ has generalized Whittaker type
$(\boldsymbol\ell;\boldsymbol\gamma)$.

\begin{proposition}\label{prop:generalized-tensor-realization}
Let $p\geq1$ and $(\ell_{2p},\gamma_{3p})\neq(0,0)$.
For each simple root $d$ of
\begin{equation}\label{gw:factor-cubic}
f(x)=4x^3-6\ell_{2p}x+3\sqrt6\gamma_{3p},
\end{equation}
there is a unique pair of parameter lists
$\boldsymbol t,\boldsymbol b$ satisfying
\eqref{gw:ell}--\eqref{gw:gamma} and $b_p=d$.
They satisfy $2t_{2p}-3b_p^2\neq0$, and
\begin{equation}\label{gw:realization}
\mathbb M^k(\boldsymbol\ell;\boldsymbol\gamma)
\cong M_p^{\Vir}(c^{\Vir},\boldsymbol t)\otimes H_p(\boldsymbol b).
\end{equation}
In particular, the universal module exists and has the basis
in Definition~\ref{def:PBW-whittaker}.
\end{proposition}
\begin{proof}
The polynomial $f$ has a simple root: otherwise it would have a
triple root, and the absence of an $x^2$ term would imply
$f(x)=4x^3$, contrary to the hypothesis.
For a simple root $d$, set $t=\ell_{2p}-d^2/2$.
The equations at the largest indices are
\begin{equation}\label{gw:top}
\ell_{2p}=t+\frac{d^2}{2},\qquad
\gamma_{3p}=\kappa dt-\frac{d^3}{3\sqrt6},\qquad
2t-3d^2=-\frac13f'(d)\neq0.
\end{equation}
Hence $t_{2p}=t$ and $b_p=d$.

For $1\leq r\leq p$, equations \eqref{gw:ell}--\eqref{gw:gamma}
have the form
\[
\binom{\ell_{2p-r}}{\gamma_{3p-r}}
=
\begin{pmatrix}
1&d\\
\kappa d&\kappa t-d^2/\sqrt6
\end{pmatrix}
\binom{t_{2p-r}}{b_{p-r}}
+\binom{E_r}{G_r},
\]
where $E_r,G_r$ depend only on $t_{2p-j},b_{p-j}$ with $j<r$.
The determinant of this matrix is $(2t-3d^2)/\sqrt6$.
Induction on $r$ therefore gives existence and uniqueness of
the remaining parameters.

Proposition~\ref{prop:general-tensor-basis} gives cyclicity
and the ordered basis of the resulting tensor product.
Its universality follows from
Proposition~\ref{prop:PBW-implies-universal}.
\end{proof}

\begin{corollary}\label{cor:generalized-tensor-quotients}
Fix the parameters in Proposition~\ref{prop:generalized-tensor-realization}.
Every generalized Whittaker module of type
$(\boldsymbol\ell;\boldsymbol\gamma)$ is isomorphic to
$
N\otimes H_p(\boldsymbol b),
$ 
where $N$ is a nonzero quotient of
$M_p^{\Vir}(c^{\Vir},\boldsymbol t)$.
This tensor product is irreducible if and only if $N$ is
irreducible.
\end{corollary}
\begin{proof}
By universality, such a module is a quotient of
\eqref{gw:realization}. Theorem~\ref{thm:tensor-submodules}
identifies the kernel with $M'\otimes H_p(\boldsymbol b)$
for a $\Vir$-submodule $M'$ of the first factor.
Taking $N=M_p^{\Vir}(c^{\Vir},\boldsymbol t)/M'$ proves
the assertion, including the irreducibility statement.
\end{proof}

\subsection{Irreducibility of universal generalized Whittaker modules}

\begin{theorem}\label{gw:criterion}
Let $p\geq1$ and $(\ell_{2p},\gamma_{3p})\neq(0,0)$, and put
\[
\Delta=\ell_{2p}^{\,3}-\frac{27}{4}\gamma_{3p}^{\,2},\qquad
\Theta=2\ell_{2p}\gamma_{3p-1}-3\gamma_{3p}\ell_{2p-1}.
\]
The universal generalized Whittaker module
$\mathbb M^k(\boldsymbol\ell;\boldsymbol\gamma)$ is
irreducible if and only if
\[
\Delta\neq0\quad\text{or}\quad\Theta\neq0.
\]
\end{theorem}
\begin{proof}
In \eqref{gw:realization}, put
$t=t_{2p}$, $s=t_{2p-1}$, $d=b_p$ and $D=2t-3d^2\neq0$.
Corollary~\ref{cor:general-tensor-simple} gives
\[
\mathbb M^k(\boldsymbol\ell;\boldsymbol\gamma)
\text{ is irreducible}\quad\Longleftrightarrow\quad
(t,s)\neq(0,0).
\]
The identities \eqref{gw:top} imply
\begin{equation}\label{gw:discriminant}
\Delta=\frac14tD^2.
\end{equation}
Thus $\Delta\neq0$ precisely when $t\neq0$.

If $\Delta=0$, then $t=0$ and $d\neq0$.
Set $e=b_{p-1}-\delta_{p,1}\sqrt6\alpha_0$.
Equations \eqref{gw:ell}--\eqref{gw:gamma} become
\[
\ell_{2p-1}=s+de,\qquad
\gamma_{3p-1}=\kappa ds-\frac{d^2e}{\sqrt6}.
\]
Together with $\ell_{2p}=d^2/2$ and
$\gamma_{3p}=-d^3/(3\sqrt6)$, these give
\[
\Theta=\frac{3d^3}{\sqrt6}s.
\]
Hence, when $\Delta=0$, the condition $\Theta\neq0$
is precisely $s\neq0$.
\end{proof}

If $\Delta\neq0$ or $\Theta\neq0$, every generalized Whittaker
module of the prescribed type is isomorphic to
$\mathbb M^k(\boldsymbol\ell;\boldsymbol\gamma)$, since it is a
nonzero quotient of this irreducible universal module.

\begin{corollary}\label{cor:rank-two-universal}
The generalized Whittaker module constructed in
Theorem~\ref{thm:rank-two} is universal and has the basis in
Definition~\ref{def:PBW-whittaker}.
\end{corollary}
\begin{proof}
By Proposition~\ref{prop:generalized-tensor-realization} and
Theorem~\ref{gw:criterion}, the universal module with these
parameters exists and is irreducible. Its homomorphism to the
module in Theorem~\ref{thm:rank-two} is surjective by cyclicity,
and hence is an isomorphism.
\end{proof}

\begin{remark}
The hypothesis $(\ell_{2p},\gamma_{3p})\neq(0,0)$ is used
in both the factorization and the criterion above.
If both parameters vanish, \eqref{gw:factor-cubic} has only
the triple root $0$. The Whittaker modules in
Theorem~\ref{thm:whittaker-universal}
belong to this case for $p=1$ and are treated separately.
\end{remark}

\subsection{Generalized Whittaker $\Walg^k(\mathfrak{sl}_3)$-modules of order one}

For the prescribed type $(\ell_1,\ell_2;\gamma_2,\gamma_3)$
with $(\ell_2,\gamma_3)\neq(0,0)$, the tensor factors can
be written explicitly. Let $d$ be a simple root of
\[
4d^3-6\ell_2d+3\sqrt6\gamma_3=0,
\]
and set
\[
\begin{aligned}
t_2&=\ell_2-\frac{d^2}{2},&
D&=2t_2-3d^2,\\
t_1&=\frac{(2t_2-d^2)\ell_1-\sqrt6d\gamma_2}{D},&
\beta&=\sqrt6\alpha_0+
       \frac{\sqrt6\gamma_2-2d\ell_1}{D}.
\end{aligned}
\]
Then
\[
\mathbb M^k(\ell_1,\ell_2;\gamma_2,\gamma_3)
\cong
M_1^{\Vir}(c^{\Vir},(t_1,t_2))\otimes H_1(\beta,d).
\]
Indeed, the equations for the two remaining parameters are
\[
\ell_1=t_1+d(\beta-\sqrt6\alpha_0),\qquad
\sqrt6\gamma_2=
2dt_1+(2t_2-d^2)(\beta-\sqrt6\alpha_0),
\]
whose determinant is $D\neq0$.

In this case, the criterion is
\[
\ell_2^3-\frac{27}{4}\gamma_3^2\neq0
\quad\text{or}\quad
2\ell_2\gamma_2-3\gamma_3\ell_1\neq0,
\]
and is precisely the condition $(t_1,t_2)\neq(0,0)$.
If both expressions vanish, the first factor is
$M_1^{\Vir}(c^{\Vir},(0,0))$, and the universal module
is reducible.

The value $d=0$ is allowed. For $\gamma_3=0$ and
$\ell_2\neq0$, this choice gives $t_2=\ell_2$,
$t_1=\ell_1$ and $H_1(\beta,0)=F_\beta$.

\subsection{Universal Whittaker $\Walg^k(\mathfrak{sl}_3)$-modules}

\begin{theorem}\label{thm:whittaker-universal}
Let $\chi_1\neq0$ and $\chi_2\in\mathbb C$, and put
$
\beta=\sqrt{3/2}\left(\frac{\chi_2}{\chi_1}
                         +\frac32\alpha_0\right)$. 
The Whittaker module
$L^{\Vir}(c^{\Vir},\chi_1)\otimes F_\beta$
constructed in Theorem~\ref{thm:vir-heis} is universal and
irreducible and has the basis in Definition~\ref{def:PBW-whittaker}.
\end{theorem}
\begin{proof}
The ordered family spans by Lemma~\ref{lem:Whittaker-spanning}
and the cyclicity proved in Theorem~\ref{thm:vir-heis}.

Put $w=u_{\chi_1}\otimes v_\beta$ and $\kappa=\sqrt{2/3}$.
For a tensor basis vector
\[
v=T(1-n_1)\cdots T(1-n_r)J(-m_1)\cdots J(-m_s)w,
\]
put $\tau(v)=(N,s,S)$, where $N=\sum_i n_i+\sum_jm_j$
and $S=\sum_jm_j$. Order these triples lexicographically,
and let $V_{<\tau}$ be the span of basis vectors with smaller
triples. Lemma~\ref{lem:vir-heis}\textup{(i)} gives
\[
\begin{aligned}
L(1-n)v&\equiv T(1-n)v
             &&\pmod{V_{<(N+n,s,S)}},\\
W(1-n)v&\equiv\kappa\chi_1J(-n)v
             &&\pmod{V_{<(N+n,s+1,S+n)}}.
\end{aligned}
\]
Indeed, the other terms of total degree $N+n$ in the second
formula have triples $(N+n,s,S)$,
$(N+n,s+1,S+j)$ with $j<n$, or $(N+n,s-1,S-q)$ with $q\geq1$.
Reordering the $T$-factors changes only terms of smaller total
degree, since
\[
[T(1-a),T(1-b)]=(b-a)T(2-a-b).
\]
It follows that the leading tensor basis vector of
\[
L(1-n_1)\cdots L(1-n_r)W(1-m_1)\cdots W(1-m_s)w
\]
is
\[
(\kappa\chi_1)^s
T(1-n_1)\cdots T(1-n_r)J(-m_1)\cdots J(-m_s)w.
\]
These leading vectors are distinct and have nonzero coefficients.
This proves linear independence.

Universality follows from Proposition~\ref{prop:PBW-implies-universal};
irreducibility was proved in Theorem~\ref{thm:vir-heis}.
\end{proof}

For $\chi_1\neq0$, every nonzero module generated by a Whittaker
vector of type $(\chi_1,\chi_2)$ is therefore isomorphic to
$L^{\Vir}(c^{\Vir},\chi_1)\otimes F_\beta$.

\begin{remark}
The case $\chi_1=0$, $\chi_2\neq0$ is usually called the
Gaiotto case. Whittaker vectors of this type were constructed
in completions of generic Verma modules in
\cite[Proposition~3.16]{BBCC24}.
Gaiotto vectors encode Nekrasov instanton partition functions
through the Kac--Shapovalov pairing \cite{SV13,BFN16}, and their
coefficients can be computed by topological recursion
\cite{BBCC24,BCU25}.

We plan to study the corresponding Whittaker modules in our forthcoming papers.
\end{remark}

\section*{Acknowledgment}

The authors acknowledge the use of OpenAI’s ChatGPT to discuss mathematical arguments, identify potential gaps and inconsistencies, and improve the language and organization of the manuscript.  All mathematical arguments and statements were independently checked by the authors, who take full responsibility for the final manuscript.

D.A. was partially supported by the Croatian Science Foundation under the project IP-2022-10-9006  and by the project “Implementation of cutting-edge research and its application as part of the Scientific Center of
Excellence for Quantum and Complex Systems, and Representations of Lie Algebras”, Grant No. PK.1.1.10.0004, co-financed by the European Union through the
European Regional Development Fund - Competitiveness and Cohesion Programme 2021-2027. Part of the research reported in this paper was conducted during D.A.’s visit to the School of Mathematical Sciences,  Xiamen University in October  2026. He is grateful to Xiamen University  for its warm hospitality. 

N.Y. was partially supported by the Natural Science Foundation of Xiamen Municipality (Grant No. 3502Z202473007), the Natural Science Foundation of Fujian Province (Grant No. 2024J01027), and the National Natural Science Foundation of China (Grant Nos. 12571031 and 12131018). Part of this research was carried out during N.Y.'s visit to the University of Zagreb in 2026. She thanks the University of Zagreb for its warm hospitality.

\end{document}